\documentclass[leqno]{article}%
\usepackage{amsfonts}
\usepackage{amsmath}
\usepackage{amsthm}
\usepackage{amssymb}
\usepackage{amsfonts}
\usepackage{xspace}
\usepackage[english,francais]{babel}%

\theoremstyle{plain}
\newtheorem{theorem}{Th\'{e}or\`{e}me}[section]
\newtheorem*{etheorem}{Theorem}
\newtheorem{lemma}[theorem]{Lemme}
\newtheorem{proposition}[theorem]{Proposition}

\newtheorem{definition}{D\'efinition}[section]

\newtheorem*{criterion}{Crit\`ere}
\newtheorem*{ecriterion}{Criterion}
\theoremstyle{remark}

\newtheorem*{notation}{Notations}

\renewenvironment{proof}{\par \noindent \textbf{D\'emonstration.}\ }{\hfill
\qedsymbol \smallskip}
\numberwithin{equation}{section}
\begin{document}

\title{Polyn\^{o}me de Hua, noyau de Bergman \\des domaines de Cartan-Hartogs \\et probl\`{e}me de Lu Qikeng}
\date{20 janvier 2007}
\author{F.~Zohra Demmad--Abdessameud \thanks{D\'{e}partement de math\'{e}matiques,
Facult\'{e} des Sciences, Universit\'{e} Saad Dahlab, Route de Soum\^{a}a, BP
270, Blida, Alg\'{e}rie ; \texttt{fz\_demmad@mail.univ-blida.dz,
fz\_demmad@yahoo.fr}}}
\maketitle

\begin{abstract}
R\'{e}duction du probl\`{e}me de Lu Qikeng pour les domaines de Cartan-Hartogs
$\widehat{\Omega}_{m}(\mu)$ \`{a} un probl\`{e}me alg\'{e}brique sur les
polyn\^{o}mes de Hua. Solution compl\`{e}te du probl\`{e}me de Lu Qikeng quand
le domaine de base $\Omega$ est un domaine sym\'{e}trique de dimension
inf\'{e}rieure ou \'{e}gale \`{a} $4$.

\emph{Classification AMS (2000)} : 32M15, 32A36. \emph{Mots cl\'{e}s} :
Domaines de Cartan, noyau de Bergman, polyn\^{o}me de Hua, conjecture de Lu Qikeng.

\end{abstract}
\tableofcontents

\selectlanguage{english}

\section*{R\'{e}sum\'{e} in English language}
\addcontentsline{toc}{section}{\protect R\'{e}sum\'{e} in English language}

The \emph{Lu Qikeng problem} for a domain $U\subset$ $\mathbb{\mathbb{C} }%
^{n}$ consists in deciding whether the Bergman kernel $K_{U}(z,w)$ of this
domain may vanish at some points of $U\times U$. A domain $U$ is called a
\emph{Lu Qikeng domain} if its Bergman kernel is zero-free on $U\times U$.

Let $\Omega$ be an irreducible bounded circled homogeneous domain and $N(z,t)$
its generic norm. For $\mu>0$ and $m$ a positive integer, let
\[
\widehat{\Omega}_{m}\left(  \mu\right)  =\left\{  \left(  z,Z\right)
\in\Omega\times\mathbb{C} ^{m},\quad\left\Vert Z\right\Vert ^{2}<N\left(
z,z\right)  ^{\mu}\right\}  .
\]
The domain $\widehat{\Omega}_{m}\left(  \mu\right)  $ is called
\emph{Cartan--Hartogs domain} (with \emph{base} $\Omega$, \emph{exponent}
$\mu$, \emph{fiber dimension} $m$).

The Bergman kernel of this domain may be explicitly computed from the generic
norm and the \emph{Hua polynomial} of $\Omega$. If $(a,b,r)$ are the numerical
invariants (multiplicities and rank) of the domain $\Omega$, its Hua
polynomial is
\[
\chi\left(  s\right)  =\chi_{a,b,r}\left(  s\right)  =\underset{j=1}%
{\overset{r}{\prod}}\left(  s+1+\left(  j-1\right)  \textstyle\frac{a}%
{2}\right)  _{1+b+\left(  r-j\right)  a},
\]
where $\left(  s+1\right)  _{k}=\prod_{i=1}^{k}(s+i)$ denotes the raising
factorial. This polynomial is related to the \emph{Hua integral} by
\[
\int_{\Omega}N(z,z)^{s}\omega(z)=\frac{\chi(0)}{\chi(s)}\int_{\Omega}%
\omega\qquad(\operatorname{Re}s>-1)
\]
(cf. \cite{YinRoos2003}). The decomposition
\[
\chi\left(  k\mu\right)  =\sum_{j=0}^{d}\mu^{j}C_{d-j}\left(  \mu\right)
\left(  k+1\right)  _{j}%
\]
of $\chi\left(  k\mu\right)  $ along raising factorials w.r.~to $k$ defines
polynomials $C_{j}(\mu)$, which are of degree $j$ in $\mu$. For $m$ a positive
integer and $\mu>0$, define%
\[
P_{\mu}^{m}(\eta)=\sum_{j=0}^{d}(m+1)_{j}C_{d-j}(\mu)\mu^{j}\eta^{j}.
\]
Note that the degree of this polynomial w.r.~to $\eta$ or $\mu$ is equal to
the dimension $d$ of $\Omega$.

The Bergman kernel of $\widehat{\Omega}_{m}\left(  \mu\right)  $ is then (cf.
\cite{YinRoos2003}, \cite{Roos2004})
\[
\widehat{K}_{m,\mu}\left(  (z,Z),(w,W)\right)  =\frac{C}{N\left(  z,w\right)
^{g+m\mu}}\eta^{m+1}P_{\mu}^{m}(\eta),
\]
where $\xi,\eta:\widehat{\Omega}_{m}\left(  \mu\right)  \times\widehat{\Omega
}_{m}\left(  \mu\right)  \rightarrow\mathbb{C}$ are defined by%
\[
\xi\left(  (z,Z),(w,W)\right)  =\frac{\left\langle Z,W\right\rangle
}{N(z,w)^{\mu}},\qquad\eta=\frac{1}{1-\xi}.
\]
\emph{The range of }$\xi$\emph{ is the unit disc }$\Delta\subset\mathbb{C}%
$\emph{ and the range of }$\eta$\emph{ is the half-plane }$\left\{
\operatorname{Re}\eta>\frac{1}{2}\right\}  $.

Thus the Lu Qikeng problem for $\widehat{\Omega}_{m}\left(  \mu\right)  $ is
reduced to the localization of the roots of $P_{\mu}^{m}$ :

\begin{etheorem}
The domain $\widehat{\Omega}_{m}\left(  \mu\right)  $ is a Lu Qikeng domain if
and only if all roots of $P_{\mu}^{m}$ are located in $\left\{
\operatorname{Re}\eta\leq\frac{1}{2}\right\}  $.
\end{etheorem}

Applying this theorem, the Lu Qikeng problem is completely solved in this
paper for all $m$, $\mu>0$ when the base is of dimension $d\leq4$. This
provides a lot of examples of Lu Qikeng and non-Lu Qikeng domains. In contrast
with the generic case of a bounded domain (\textquotedblleft The Lu Qikeng
conjecture fails generically\textquotedblright, see \cite{Boas1996}),
\textquotedblleft most\textquotedblright\ of the domains $\widehat{\Omega}%
_{m}\left(  \mu\right)  $ are Lu Qikeng domains. Actually, the domain
$\widehat{\Omega}_{m}\left(  \mu\right)  $ is a Lu Qikeng domain for $m\geq
m_{\Omega}$ and for all $\mu>0$, where $m_{\Omega}$ is an integer depending on
the base $\Omega$; for $1\leq m<m_{\Omega}$, there exists a positive real
number $\mu_{m}$ such that the domain $\widehat{\Omega}_{m}\left(  \mu\right)
$ is a Lu Qikeng domain if and only $0<\mu\leq\mu_{m}$.

\bigskip

Results are as follows: \smallskip

1. If $\Omega$ is the unit disc $\Delta\subset\mathbb{C}$, the domain
$\widehat{\Omega}_{m}\left(  \mu\right)  $ is a Lu Qikeng domain for all
$m\geq1$ and all $\mu>0$. This is recalled here only for sake of completeness.

\smallskip

2. If $\Omega$ is the Hermitian ball of dimension $2$, the domain
$\widehat{\Omega}_{m}\left(  \mu\right)  $ is a Lu Qikeng domain if and only if

\begin{enumerate}
\item $m=1$, $\mu\leq2$;

\item $m=2$, $\mu\leq4$;

\item $m\geq m_{\Omega}=3$, for all $\mu>0$.
\end{enumerate}

For $m=1$, the result is due to H.P. Boas, Siqi Fu, E. Straube
\cite{BoasFuStraube1999}; see also \cite{Yin2006}. For $m>1$, results are new.

\smallskip

3. If $\Omega$ is the Hermitian ball of dimension $3$, for $1\leq m\leq5$, the
polynomial $q_{m}$ (of degree $2$ or $3$) defined by
\[
q_{m}(\mu)=P_{\mu}^{m}\left(  \textstyle\frac{1}{2}\right)
\]
has a unique positive root $\mu_{m}$ and%
\[
0<\mu_{1}=\sqrt{2}<\mu_{2}<\mu_{3}<\mu_{4}<\mu_{5}.
\]
The domain $\widehat{\Omega}_{m}\left(  \mu\right)  $ is a Lu Qikeng domain if
and only if

\begin{enumerate}
\item $1\leq m\leq5$, $0<\mu\leq\mu_{m}$;

\item $m\geq m_{\Omega}=6$, for all $\mu>0$.
\end{enumerate}

For $m=1$, this result has been obtained by Weiping Yin \cite{Yin2006}, by a
slightly different method. Results are new for $m>1$.

\smallskip

4. If $\Omega$ is the Lie ball of dimension $3$ (domain of type $IV_{3}\simeq
III_{2}$, symmetric matrices), the same type of result holds as in the
preceding case, with different $q_{m}$ and $\mu_{m}$,
\[
0<\mu_{1}=\textstyle\frac{2}{\sqrt{3}}<\mu_{2}<\mu_{3}<\mu_{4}<\mu_{5},
\]
but again $m_{\Omega}=6$. Here the results are new for all $m$.

\smallskip

5. If $\Omega$ is the Hermitian ball of dimension $4$, the polynomial $q_{m}$
defined by
\[
q_{m}(\mu)=P_{\mu}^{m}\left(  \textstyle\frac{1}{2}\right)
\]

\begin{itemize}
\item has two positive roots $\mu_{m}=\mu_{m,1}<\mu_{m,2}$ for $m=1,2$;

\item has one positive root $\mu_{m}$ for $3\leq m\leq7$;

\item is positive for all $\mu>0$ if $m\geq m_{\Omega}=8$.
\end{itemize}

Moreover,
\[
0<\mu_{1}=\textstyle\sqrt{\frac{3}{2}}<\mu_{2}<\mu_{3}<\mu_{4}<\mu_{5}<\mu
_{6}<\mu_{7}.
\]
The domain $\widehat{\Omega}_{m}\left(  \mu\right)  $ is a Lu Qikeng domain if
and only if

\begin{enumerate}
\item $1\leq m\leq7$, $0<\mu\leq\mu_{m}$;

\item $m\geq m_{\Omega}=8$, for all $\mu>0$.
\end{enumerate}

For $m=1$, this result has been obtained by Jong-do Park and Liyou Zhang
(2006, unpublished). Results are new for $m>1$.

\smallskip

6. If $\Omega$ is the Lie ball of dimension $4$ (domain of type $IV_{4}\simeq
I_{2,2}$, $2\times2$ matrices), the same type of result as in the preceding
case holds, with different $q_{m}$ and $\mu_{m}$,
\[
0<\mu_{1}=\textstyle\frac{1}{2}\sqrt{23-\sqrt{337}}<\mu_{2}<\mu_{3}<\mu
_{4}<\mu_{5}<\mu_{6}<\mu_{7},
\]
and $m_{\Omega}=8$. Here the results are new for all $m$.

\bigskip Results may be summarized in the following theorem:

\begin{etheorem}
Let $\Omega$ be an irreducible bounded circled homogeneous domain of dimension
at most $4$. Then the polynomial
\[
q_{m}(\mu)=P_{\mu}^{m}\left(  \textstyle\frac{1}{2}\right)
\]
has $0$, $1$ or $2$ positive roots. If $q_{m}$ has no positive root, let
$\mu_{m}=+\infty$; if $q_{m}$ has one positive root, denote this root by
$\mu_{m}=\mu_{m,1}$ and let $\mu_{m,2}=+\infty$; if $q_{m}$has two positive
roots, denote these roots by $\mu_{m,1}$, $\mu_{m,2}$ and let $\mu_{m}%
=\mu_{m,1}<\mu_{m,2}$.

The polynomial $P_{\mu}^{m}$ has

\begin{itemize}
\item no root in $\left\{  \operatorname{Re}\eta>\frac{1}{2}\right\}  $ if
$0<\mu\leq\mu_{m}$;

\item one root in $\left\{  \operatorname{Re}\eta>\frac{1}{2}\right\}  $ if
$\mu_{m,1}<\mu\leq\mu_{m,2}$;

\item two roots in $\left\{  \operatorname{Re}\eta>\frac{1}{2}\right\}  $ if
$\mu_{m,2}<\mu$.
\end{itemize}

The Cartan--Hartogs domain $\widehat{\Omega}_{m}\left(  \mu\right)  $ is a Lu
Qikeng domain if and only if%
\[
0<\mu\leq\mu_{m}.
\]

\end{etheorem}

The values of the positive roots of $q_{m}$ are given in the following table.

\bigskip%

\begin{tabular}
[c]{|c|c|c|c|c|c|}\hline
Type & $I_{1,2}$ & $I_{1,3}\simeq II_{3}$ & $III_{2}\simeq IV_{3}$ & $I_{1,4}$
& $I_{2,2}\simeq IV_{4}$\\\hline
$\mu_{1,1}$ & $2$ & $\sqrt{2}$ & $\frac{2}{\sqrt{3}}$ & $\sqrt{\frac{3}{2}}$ &
$\simeq1.07732$\\\hline
$\mu_{1,2}$ & $+\infty$ & $+\infty$ & $+\infty$ & $4$ & $\simeq3.21549$%
\\\hline
$\mu_{2,1}$ & $4$ & $\frac{1+\sqrt{7}}{2}$ & $\frac{3+\sqrt{73}}{8}$ &
$\simeq1.41518$ & $\simeq1.21176$\\\hline
$\mu_{2,2}$ & $+\infty$ & $+\infty$ & $+\infty$ & $\simeq11.333$ &
$\simeq9.08062$\\\hline
$\mu_{3}$ & $+\infty$ & $1+\sqrt{\frac{5}{2}}$ & $2$ & $\simeq1.61819$ &
$\simeq1.41824$\\\hline
$\mu_{4}$ & $+\infty$ & $2+\sqrt{6}$ & $\frac{1}{2}\left(  3+\sqrt{\frac
{43}{3}}\right)  $ & $\simeq2.10335$ & $\simeq1.74173$\\\hline
$\mu_{5}$ & $+\infty$ & $8+\sqrt{70}$ & $2\left(  3+\sqrt{10}\right)  $ &
$\simeq2.8029$ & $\simeq2.29476$\\\hline
$\mu_{6}$ & $+\infty$ & $+\infty$ & $+\infty$ & $\simeq4.22107$ &
$\simeq3.42405$\\\hline
$\mu_{7}$ & $+\infty$ & $+\infty$ & $+\infty$ & $\simeq8.60867$ &
$\simeq6.92986$\\\hline
\end{tabular}

\bigskip

Proofs are by case-by-case computation and study of the polynomial $P_{\mu
}^{m}$. For base domain $\Omega$ of dimension $3$ or $4$, most computations
have been done with \textsc{Mathematica}; these computations involve only
algebraic operations on polynomials and localization of their roots, and they
can be checked with any software for symbolic calculus. The localization of
the roots of $P_{\mu}^{m}$ w.r.~to $\left\{  \operatorname{Re}\eta=\frac{1}%
{2}\right\}  $ is easy in degree $1$ or $2$. In degrees $3$ and $4$, we use
the following criteria, which result from the \emph{Routh--Hurwitz criterion}
(see \cite[Chap. 15]{Gantmacher1965} and Propositions \ref{LocDegre3},
\ref{LocDegre4}).

\begin{ecriterion}
[Degree $3$]Let $P(z)=\alpha+\beta z+\gamma z^{2}+\delta z^{3}$ be a
polynomial with real coefficients and $\delta>0$. Then all roots of $P$ are
located in $\left\{  \operatorname{Re}z<\frac{1}{2}\right\}  $ if and only if
\[
P\left(  \textstyle\frac{1}{2}\right)  >0,\quad P^{\prime}\left(
\textstyle\frac{1}{2}\right)  >0,\quad\Delta_{2}\equiv\left(  \gamma
+\delta\right)  \left(  \beta+\gamma+\delta\right)  -\alpha\delta>0.
\]

\end{ecriterion}

\begin{ecriterion}
[Degree $4$]Let
\[
P(z)=\alpha+\beta z+\gamma z^{2}+\delta z^{3}+\varepsilon z^{4}%
\]
be a polynomial with real coefficients and $\varepsilon>0$. Then all roots of
$P$ are located in $\left\{  \operatorname{Re}z<\frac{1}{2}\right\}  $ if and
only if%
\begin{align*}
P  &  \left(  \textstyle\frac{1}{2}\right)  >0,\quad P^{\prime}\left(
\textstyle\frac{1}{2}\right)  >0,\quad P^{\prime\prime}\left(  \textstyle\frac
{1}{2}\right)  >0,\\
\Delta_{3}  &  \equiv\left(  \varepsilon+\delta+\gamma+\beta\right)  \left[
\left(  \varepsilon+\delta+\gamma\right)  \left(  \varepsilon+\delta\right)
-\varepsilon\beta\right]  -\left(  2\varepsilon+\delta\right)  ^{2}\alpha>0.
\end{align*}

\end{ecriterion}

\emph{Sketch of proofs.}

\smallskip

1. If $\Omega$ is the unit disc, $P_{\mu}^{m}$ has degree $1$ and
\[
q_{m}(\mu)=P_{\mu}^{m}\left(  \textstyle\frac{1}{2}\right)  =\textstyle\frac
{(m-1)\mu}{2}+1
\]
is always positive.

\smallskip

2. If $\Omega$ is the Hermitian ball of dimension $2$, the polynomial $P_{\mu
}^{m}(\eta)$ has degree $2$; its roots (real or imaginary conjugate) lie in
$\left\{  \operatorname{Re}\eta\leq\frac{1}{2}\right\}  $ if and only if
\[
P_{\mu}^{m}\left(  \textstyle\frac{1}{2}\right)  \geq0,\qquad\textstyle\frac
{\operatorname*{d}}{\operatorname*{d}\eta}P_{\mu}^{m}\left(  \textstyle\frac
{1}{2}\right)  \geq0.
\]
Computation shows that
\begin{align*}
q_{m}(\mu)  &  =P_{\mu}^{m}\left(  \textstyle\frac{1}{2}\right)
=\textstyle\frac{m(m-3)}{4}\mu^{2}+\textstyle\frac{3(m-1)}{2}\mu+2,\\
q_{m}^{1}(\mu)  &  =\textstyle\frac{1}{(m+1)\mu}\frac{\operatorname*{d}%
}{\operatorname*{d}\eta}P_{\mu}^{m}\left(  \textstyle\frac{1}{2}\right)
=3+(m-1)\mu;
\end{align*}
the results easily follows by inspection of the special cases $m=1$ and $m=2$.

\smallskip

3. If $\Omega$ is the Hermitian ball of dimension $3$, the polynomial $P_{\mu
}^{m}(\eta)$ has degree $3$. According to the above criterion for degree $3$,
all roots of $P_{\mu}^{m}(\eta)=\alpha+\beta\eta+\gamma\eta^{2}+\delta\eta
^{3}$ are located in $\left\{  \operatorname{Re}\eta<\frac{1}{2}\right\}  $ if
and only if
\begin{align*}
P_{\mu}^{m}  &  \left(  \textstyle\frac{1}{2}\right)  >0,\quad\frac
{\operatorname*{d}}{\operatorname*{d}\eta}P_{\mu}^{m}\left(  \textstyle\frac
{1}{2}\right)  >0,\\
\Delta_{2}  &  =\left(  \gamma+\delta\right)  \left(  \beta+\gamma
+\delta\right)  -\alpha\delta>0.
\end{align*}
Direct computations show that the first condition implies the second. A study
of the polynomials
\[
q_{m}(\mu)=P_{\mu}^{m}\left(  \textstyle\frac{1}{2}\right)  =6+\textstyle\frac
{11\left(  m-1\right)  }{2}\mu+\frac{3m(m-3)}{2}\mu^{2}+\frac{\left(
m-1\right)  \left(  m^{2}-5m-2\right)  }{8}\mu^{3}%
\]
and their comparison give the results about the $\mu_{m}$'s. Computation of
the third condition gives rise to a polynomial condition on $\mu$ and $m$,
which again appears to be implied by
\[
q_{m}(\mu)\geq0.
\]

\smallskip

4. If $\Omega$ is the Lie ball of dimension $3$, the treatment is entirely
analogous to the case of type $I_{1,3}$, with
\[
q_{m}(\mu)=\textstyle\frac{1}{8}\left(  m-1\right)  \left(  m^{2}-5m-2\right)
\mu^{3}+\frac{9}{8}m\left(  m-3\right)  \mu^{2}+\frac{13}{4}\left(
m-1\right)  \mu+3.
\]

\smallskip

5. If $\Omega$ is the Hermitian ball of dimension $4$, the polynomial $P_{\mu
}^{m}(\eta)$ has degree $4$. According to the above criterion for degree $4$,
all roots of $P_{\mu}^{m}(\eta)=\alpha+\beta\eta+\gamma\eta^{2}+\delta\eta
^{3}+\varepsilon\eta^{4}$ are located in $\left\{  \operatorname{Re}\eta
<\frac{1}{2}\right\}  $ if and only if
\begin{align*}
P_{\mu}^{m}  &  \left(  \textstyle\frac{1}{2}\right)  >0,\quad\textstyle\frac
{\operatorname*{d}}{\operatorname*{d}\eta}P_{\mu}^{m}\left(  \textstyle\frac
{1}{2}\right)  >0,\quad\textstyle\frac{\operatorname*{d}^{2}}%
{\operatorname*{d}\eta^{2}}P_{\mu}^{m}\left(  \textstyle\frac{1}{2}\right)
>0,\\
\Delta_{3}  &  \equiv\left(  \varepsilon+\delta+\gamma+\beta\right)  \left[
\left(  \varepsilon+\delta+\gamma\right)  \left(  \varepsilon+\delta\right)
-\varepsilon\beta\right]  -\left(  2\varepsilon+\delta\right)  ^{2}\alpha>0.
\end{align*}
The polynomial $q_{m}(\mu)=P_{\mu}^{m}\left(  \frac{1}{2}\right)  $ is equal
to
\begin{align*}
q_{m}(\mu)  &  =P_{\mu}^{m}\left(  \textstyle\frac{1}{2}\right) \\
&  =(4+m\mu)\left(  6+\textstyle\frac{19m-25}{4}\mu+m(m-5)\mu^{2}+\frac
{m^{3}-10m^{2}+15m+10}{16}\mu^{3}\right)  .
\end{align*}
This polynomial has

\begin{itemize}
\item two positive roots $\mu_{m}=\mu_{m,1}<\mu_{m,2}$ for $m=1,2$;

\item one positive root $\mu_{m}$ for $3\leq m\leq7$;

\item no positive root for $m\geq8$.
\end{itemize}

A case-by-case study shows that $0<\mu\leq\mu_{m}$ implies
\[
\textstyle\frac{\operatorname*{d}}{\operatorname*{d}\eta}P_{\mu}^{m}\left(
\frac{1}{2}\right)  >0,\quad\frac{\operatorname*{d}^{2}}{\operatorname*{d}%
\eta^{2}}P_{\mu}^{m}\left(  \frac{1}{2}\right)  >0.
\]
Also, the condition $\Delta_{3}>0$ applied to $P_{\mu}^{m}$ gives a polynomial
condition (in general of degree $6$ w.r.~to $\mu$), which is shown to be
satisfied for all $\mu$ such that $\frac{\operatorname*{d}}{\operatorname*{d}%
\eta}P_{\mu}^{m}\left(  \frac{1}{2}\right)  \geq0$. Finally, $P_{\mu}^{m}$ has
no root in $\left\{  \operatorname{Re}\eta>\frac{1}{2}\right\}  $ if and only
if $0<\mu\leq\mu_{m}$ (with $\mu_{m}=+\infty$ for $m\geq m_{\Omega}=8$). It is
also possible to check that the number of roots of $P_{\mu}^{m}$ in $\left\{
\operatorname{Re}\eta>\frac{1}{2}\right\}  $ is $1$ for $\mu_{m}<\mu\leq
\mu_{m,2}$, and $2$ for $\mu>\mu_{m,2}$ (with $\mu_{m,2}=+\infty$ for $m\geq
3$); moreover, these roots are always real.

\smallskip

6. If $\Omega$ is the Lie ball of dimension $4$, the treatment is entirely
analogous to the case of type $I_{1,4}$, with
\begin{align*}
q_{m}(\mu)=P_{\mu}^{m}\left(  \textstyle\frac{1}{2}\right)  =12  &
+14(m-1)\mu+\textstyle\frac{23m(m-3)}{4}\mu^{2}\\
&  +(m-1)(m^{2}-5m-2)\mu^{3}+\textstyle\frac{m^{3}-10m^{2}+15m+10}{16}\mu^{4}.
\end{align*}
\selectlanguage{french}

\section*{Introduction}
\addcontentsline{toc}{section}{\protect Introduction}

Le \emph{probl\`{e}me de Lu Qikeng} pour un ouvert $U$ de $\mathbb{\mathbb{C}
}^{n}$ consiste \`{a} d\'{e}terminer si le noyau de Bergman $K_{U}(z,w)$ de ce
domaine peut avoir des z\'{e}ros. Ce probl\`{e}me a \'{e}t\'{e} pos\'{e} par
Lu Qikeng en 1966. Le nom de \emph{conjecture de Lu Qikeng} a \'{e}t\'{e}
donn\'{e} (par M.~Skwarsczynski en 1969) \`{a} l'hypoth\`{e}se suivant
laquelle le noyau de Bergman d'un ouvert n'aurait pas de z\'{e}ros. Un domaine
$U$ sera appel\'{e} \emph{domaine de Lu Qikeng} si son noyau de Bergman ne
s'annule pas dans $U\times U$.

Soit $\Omega$ un domaine homog\`{e}ne born\'{e} cercl\'{e} irr\'{e}ductible,
$N(z,t)$ sa norme g\'{e}\-n\'{e}\-rique, $\chi$ son polyn\^{o}me de Hua. Pour
$\mu>0$ et $m$ entier positif, on consid\`{e}re le \emph{domaine de
Cartan--Hartogs} $\widehat{\Omega}_{m}\left(  \mu\right)  $ construit
au-dessus de $\Omega$:
\[
\widehat{\Omega}_{m}\left(  \mu\right)  =\left\{  \left(  z,Z\right)
\in\Omega\times\mathbb{C}^{m},\quad\left\Vert Z\right\Vert ^{2}<N\left(
z,z\right)  ^{\mu}\right\}  .
\]
Les domaines de Cartan--Hartogs ont \'{e}t\'{e} introduits en 1998 par G.~Roos
et Weiping Yin ; ils g\'{e}n\'{e}ralisent les ellispso\"{\i}des complexes, qui
correspondent au cas o\`{u} $\Omega$ est le disque unit\'{e} de $\mathbb{C}$.
\emph{Ces domaines sont en g\'{e}n\'{e}ral non homog\`{e}nes, mais les orbites
du groupe d'automorphismes sont alors param\'{e}tr\'{e}es par }$[0,1[$. Le
noyau de Bergman de ces domaines a \'{e}t\'{e} obtenu dans le cas
g\'{e}n\'{e}ral dans \cite{YinRoos2003} ; cf. \'{e}galement \cite{Roos2004}.

Dans cet article, nous \'{e}tudions le probl\`{e}me de Lu Qikeng pour les
domaines de Cartan--Hartogs. Nous montrons (th\'{e}or\`{e}me \ref{SolAlgPbLu})
qu'il se r\'{e}duit \`{a} la localisation, par rapport \`{a} la droite
$\left\{  \operatorname{Re}\eta=\frac{1}{2}\right\}  $, des racines d'un
polyn\^{o}me $P_{\mu}^{m}$ ; ce polyn\^{o}me, de degr\'{e} \'{e}gal \`{a} la
dimension $d$ de $\Omega$, se d\'{e}duit du polyn\^{o}me de Hua de $\Omega$
par une transformation combinatoire.

Nous appliquons ensuite ce th\'{e}or\`{e}me \`{a} la solution compl\`{e}te du
probl\`{e}me de Lu Qikeng pour les domaines $\widehat{\Omega}_{m}\left(
\mu\right)  $ lorsque $\Omega$ est un domaine sym\'{e}trique irr\'{e}ductible
de dimension au plus $4$. Les r\'{e}sultats font appara\^{\i}tre la situation
suivante, dont on conjecture qu'elle se g\'{e}n\'{e}ralise pour toute base
$\Omega$ : pour $\Omega$ et $m\geq1$ fix\'{e}s, il existe $\mu_{\Omega,m}$,
$0<$ $\mu_{\Omega,m}\leq\infty$ tel que $\widehat{\Omega}_{m}\left(
\mu\right)  $ est un domaine de Lu Qikeng si et seulement si $0<\mu\leq
\mu_{\Omega,m}$. La borne $\mu_{\Omega,m}$ est caract\'{e}ris\'{e}e comme la
plus petite racine positive du polyn\^{o}me $q_{m}(\mu)=P_{\mu}^{m}\left(
\frac{1}{2}\right)  $ ; on a $\mu_{\Omega,m}=+\infty$ pour $m$ assez grand. De
plus, si le domaine $\Omega$ n'est pas un domaine de Lu Qikeng, i.e. si
$\mu>\mu_{\Omega,m}$, il est possible de pr\'{e}ciser le nombre de racines de
$P_{\mu}^{m}$ dans $\left\{  \operatorname{Re}\eta>\frac{1}{2}\right\}  $, de
v\'{e}rifier que celles-ci sont toujours r\'{e}elles et de d\'{e}crire la
vari\'{e}t\'{e} des points de $\widehat{\Omega}_{m}\left(  \mu\right)
\times\widehat{\Omega}_{m}\left(  \mu\right)  $ o\`{u} le noyau de Bergman de
$\widehat{\Omega}_{m}\left(  \mu\right)  $ s'annule.

Les r\'{e}sultats obtenus lorsque la base $\Omega$ est un domaine
sym\'{e}trique irr\'{e}\-duc\-tible de dimension au plus $4$ fournissent un
grand nombre d'exemples de domaines de Lu Qikeng et de domaines qui n'ont pas
cette propri\'{e}t\'{e}. Contrairement au cas g\'{e}n\'{e}rique d'un domaine
born\'{e} de $\mathbb{C}^{n}$, qui n'est pas de Lu Qikeng (cf. \cite{Boas1996}%
), \og la plupart\fg des domaines $\widehat{\Omega}_{m}\left(  \mu\right)  $
sont des domaines de Lu Qikeng. En effet, le domaine $\widehat{\Omega}%
_{m}\left(  \mu\right)  $ est un domaine de Lu Qikeng pour $m\geq m_{\Omega}$
et pour tout $\mu>0$, o\`{u} $m_{\Omega}$ est un entier qui d\'{e}pend de la
base $\Omega$ ; pour $1\leq m<m_{\Omega}$, le domaine $\widehat{\Omega}%
_{m}\left(  \mu\right)  $ est de Lu Qikeng si et seulement si $0<\mu\leq
\mu_{\Omega,m}$.

\bigskip

Ce travail est organis\'{e} comme suit: Dans les sections \ref{PolynHua} et
\ref{Bergman}, nous rappelons la d\'{e}finition du polyn\^{o}me de Hua d'un
domaine sym\'{e}trique $\Omega$ et le calcul du noyau de Bergman de
$\widehat{\Omega}_{m}\left(  \mu\right)  $. La section \ref{LuProblem} est
essentiellement consacr\'{e}e \`{a} la d\'{e}monstration du th\'{e}or\`{e}me
de r\'{e}duction \ref{SolAlgPbLu}.

La section \ref{LowDim} d\'{e}crit la solution compl\`{e}te du probl\`{e}me de
Lu Qikeng pour $\widehat{\Omega}_{m}\left(  \mu\right)  $ lorsque la base
$\Omega$ est de dimension au plus $4$ ; les r\'{e}sultats g\'{e}n\'{e}ralisent
des r\'{e}sultats obtenus par Yin Weiping \cite{Yin2006}, Zhang Liyou et Park
Jong-do (2006, non publi\'{e}) lorsque $m=1$ et que $\Omega$ est une boule
hermitienne de dimension $3$ ou $4$. Pour all\'{e}ger cette section, les
r\'{e}sultats auxiliaires utilis\'{e}s dans ces cas particuliers ont
\'{e}t\'{e} regroup\'{e}s dans l'annexe \ref{Tables}. \`{A} partir de la
dimension $3$, les calculs ont \'{e}t\'{e} faits \`{a} l'aide de
\textsc{Mathematica} ; s'agissant uniquement d'op\'{e}rations alg\'{e}briques
sur les polyn\^{o}mes et de localisation de leurs racines, ces calculs peuvent
\^{e}tre v\'{e}rifi\'{e}s avec tout autre logiciel de calcul symbolique.

Enfin, l'annexe \ref{Localisation} regroupe les crit\`{e}res utilis\'{e}s pour
la localisation des racines de polyn\^{o}mes par rapport \`{a} $\left\{
\operatorname{Re}z<\frac{1}{2}\right\}  $. En degr\'{e}s $3$ et $4$, ces
crit\`{e}res sont d\'{e}duits du crit\`{e}re de stabilit\'{e} classique de
Routh--Hurwitz et du crit\`{e}re de Li\'{e}nard et Chipart.

\section{Polyn\^{o}mes du type de Hua\label{PolynHua}}

\subsection{D\'{e}finition}

On consid\`{e}re un triplet $(a,b,r)$ d'entiers naturels, avec $r>0$. Le
\emph{polyn\^{o}me du type de Hua} $\chi=\chi_{a,b,r}$ est le polyn\^{o}me
d\'{e}fini par
\begin{equation}
\chi\left(  s\right)  =\chi_{a,b,r}\left(  s\right)  =\underset{j=1}%
{\overset{r}{\prod}}\left(  s+1+\left(  j-1\right)  \textstyle\frac{a}%
{2}\right)  _{1+b+\left(  r-j\right)  a}, \label{HuaPol}%
\end{equation}

On a
\begin{equation}
\deg\chi=d=r+\textstyle\frac{r\left(  r-1\right)  }{2}a+rb. \label{dimension}%
\end{equation}
En effet, on a
\[
\deg\chi=\underset{j=1}{\overset{r}{\sum}}\left(  1+b+\left(  r-j\right)
a\right)  =r+rb+\textstyle\frac{r\left(  r-1\right)  }{2}a.
\]

Le polyn\^{o}me $\chi$ est li\'{e} \`{a} l'\emph{int\'{e}grale de Selberg}
\cite{Selberg1944}, pour $\operatorname{Re}s>-1$, par
\begin{equation}
\int_{0}^{1}\cdots\int_{0}^{1}\prod_{j=1}^{r}(1-t_{j})^{s}t_{j}^{b}%
\prod_{1\leq j<k\leq r}\left\vert t_{j}-t_{k}\right\vert ^{a}\;\mathrm{d}%
t_{1}\ldots\mathrm{d}t_{r}=\textstyle\frac{C(a,b,r)}{\chi(s)},
\label{IntDefHua}%
\end{equation}
o\`{u}
\[
C(a,b,r)=\prod_{j=1}^{r}\textstyle\frac{\Gamma(b+1+(j-1)\frac{a}{2}%
)\Gamma(j\frac{a}{2}+1)}{\Gamma(\frac{a}{2}+1)}.
\]

\subsection{Polyn\^{o}mes de Hua des domaines hermitiens sym\'{e}\-tri\-ques}

Soit $\Omega$ un domaine hermitien sym\'{e}trique irr\'{e}ductible. On
d\'{e}signe par $N$ la $\emph{norme}$ g\'{e}n\'{e}rique de $\Omega$ et par
$a,b,r$ ses invariants num\'{e}riques. L'\emph{int\'{e}grale de Hua }%
$\int_{\Omega}N(z,z)^{s}\omega(z)$ de $\Omega$ est donn\'{e}e par
\begin{equation}
\int_{\Omega}N(z,z)^{s}\omega(z)={\textstyle\frac{\chi(0)}{\chi(s)}}%
\int_{\Omega}\omega\qquad(\operatorname{Re}s>-1) \label{HuaIntegral}%
\end{equation}
(cf. \cite{YinRoos2003}), o\`{u} $\chi=\chi_{a,b,r}$.

Il est bien connu que le noyau reproduisant de l'espace \`{a} poids $H\left(
\Omega,N(z,z)^{s}\right)  $ est
\[
K_{\Omega,s}(z,t)=C_{s}N(z,t)^{-g-s},
\]
o\`{u} $C_{s}$ est une constante d\'{e}pendant de $s$. Nous rappelons
ci-dessous la d\'{e}mons\-tra\-tion et pr\'{e}cisons la relation entre $C_{s}$
et $\chi(s)$.

\begin{lemma}
\label{Bergman-poids}Soit $\Omega$ un domaine homog\`{e}ne born\'{e}
cercl\'{e} irr\'{e}ductible, $N(z,t)$ sa norme g\'{e}n\'{e}rique, $\chi$ son
polyn\^{o}me de Hua. Pour $s>0$, le noyau reproduisant de l'espace \`{a} poids
$H\left(  \Omega,N(z,z)^{s}\right)  $ est
\begin{equation}
K_{\Omega,s}(z,t)=C_{\Omega}N(z,t)^{-g-s}\chi(s), \label{Bergman-poids2}%
\end{equation}
o\`{u} $C_{\Omega}=\frac{1}{\chi\left(  0\right)  \operatorname{vol}\left(
\Omega\right)  }.$
\end{lemma}

\begin{proof}
Si $\phi$ est un automorphisme de $\Omega$ on a%
\begin{equation}
N\left(  \phi\left(  t\right)  ,\phi\left(  t\right)  \right)  ^{g}=\left\vert
J\phi\left(  t\right)  \right\vert ^{2}N\left(  t,t\right)  ^{g}.
\label{Hua-eq1}%
\end{equation}
Ceci r\'{e}sulte des relations
\begin{align*}
B\left(  \phi\left(  z\right)  ,\phi\left(  t\right)  \right)   &
=\operatorname*{d}\phi\left(  z\right)  B\left(  z,t\right)  \operatorname*{d}%
\phi^{\ast}\left(  z\right)  ,\\
\det B\left(  z,z\right)   &  =N\left(  z,z\right)  ^{g}.
\end{align*}
Soit
\[
\left\Vert f\right\Vert _{s}^{2}=\int_{\Omega}\left\vert f\left(  t\right)
\right\vert ^{2}N\left(  t,t\right)  ^{s}\omega\left(  t\right)
\]
la norme de $H\left(  \Omega,N(z,z)^{s}\right)  $. Si $\phi$ est un
automorphisme de $\Omega$, on a par changement de variable dans
l'int\'{e}grale et en appliquant (\ref{Hua-eq1}),
\begin{align*}
\left\Vert f\right\Vert _{s}^{2}  &  =\int_{\Omega}\left\vert f\circ
\phi\right\vert ^{2}N\left(  \phi\left(  t\right)  ,\phi\left(  t\right)
\right)  ^{s}\omega\left(  \phi\left(  t\right)  \right) \\
&  =\int_{\Omega}\left\vert f\circ\phi\right\vert ^{2}N\left(  t,t\right)
^{s}\left\vert J\phi\left(  t\right)  \right\vert ^{\tfrac{2s}{g}}\left\vert
J\phi\left(  t\right)  \right\vert ^{2}\omega\left(  t\right)  =\left\Vert
(f\circ\phi)\,\left(  J\phi\right)  ^{\tfrac{s}{g}+1}\right\Vert _{s}^{2}.
\end{align*}
L'application $f\longmapsto(f\circ\phi)\,\left(  J\phi\right)  ^{\tfrac{s}%
{g}+1}$ est donc un automorphisme de l'espace de Hilbert $H\left(
\Omega,N(z,z)^{s}\right)  $. Le noyau reproduisant de cet espace v\'{e}rifie
donc la relation de transformation
\[
K_{\Omega,s}(z,z)=\left\vert J\phi\left(  z\right)  \right\vert ^{\tfrac
{2s}{g}+2}K_{\Omega,s}(\phi\left(  z\right)  ,\phi\left(  z\right)  ).
\]
Soit $z\in\Omega$ et $\phi\in\operatorname*{Aut}\Omega$ tel que $\phi\left(
z\right)  =0$ ; comme $N\left(  0,0\right)  =1$ , on d\'{e}duit de
(\ref{Hua-eq1})
\[
N\left(  z,z\right)  ^{g}=\left\vert J\phi\left(  z\right)  \right\vert
^{-2},
\]
d'o\`{u} $K_{\Omega,s}(z,z)=C_{s}N(z,z)^{-g-s}$, avec $C_{s}=K_{\Omega
,s}(0,0)$. On a donc
\[
K_{\Omega,s}(z,t)=C_{s}N(z,t)^{-g-s},
\]
les deux membres \'{e}tant analytique-r\'{e}els. En particulier, $K_{\Omega
,s}(0,t)=C_{s}$ et $1=C_{s}\int N(t,t)^{s}\omega$ ; d'o\`{u}, en utilisant
(\ref{HuaIntegral}), $C_{s}=\frac{\chi(s)}{\chi(0)\operatorname{vol}\Omega
}=C_{\Omega}\chi(s)$.
\end{proof}

\section{Noyau de Bergman des domaines de Cartan-Hartogs\label{Bergman}}

\subsection{Noyau de Bergman virtuel}

Soit $V$ un espace vectoriel hermitien de dimension finie $n$, dont on note
$\left\Vert ~\right\Vert =\left\Vert ~\right\Vert _{V}$ la norme hermitienne
et $\omega_{V}(z)=\left(  \frac{\operatorname*{i}}{2\pi}\partial
\overline{\partial}\left\Vert z\right\Vert ^{2}\right)  ^{n}$ la forme volume
associ\'{e}e. Soit $\Omega$ un domaine de $V$ et $p:\Omega\rightarrow
]0,+\infty\lbrack$ une fonction continue positive sur $\Omega$. L'espace des
fonctions holomorphes sur $\Omega$ est not\'{e} $\operatorname{Hol}(\Omega)$.
On note $H(\Omega)$ l'espace de Bergman {%
\[
H(\Omega)=H\left(  \Omega,\omega_{V}\right)  =\left\{  f\in\operatorname{Hol}%
(\Omega)\mid\left\Vert f\right\Vert _{\Omega}^{2}=\int_{\Omega}\left\vert
f(z)\right\vert ^{2}\omega_{V}(z)<\infty\right\}
\]
et }$H(\Omega,p)=H(\Omega,p\omega_{V})$ l'espace de Bergman \`{a} poids{%
\[
H(\Omega,p\omega_{V})=\left\{  f\in\operatorname{Hol}(\Omega)\mid\left\Vert
f\right\Vert _{\Omega,p}^{2}=\int_{\Omega}\left\vert f(z)\right\vert
^{2}p(z)\omega_{V}(z)<\infty\right\}  .
\]
Les produits scalaires de ces espaces de Hilbert sont not\'{e}s respectivement
}$\left(  ~\mid~\right)  _{\Omega}$ et $\left(  ~\mid~\right)  _{\Omega,p}$.
Le noyau de Bergman de $\Omega$ (noyau reproduisant de $H(\Omega)$) est
not\'{e} $K_{\Omega}(z,t)$; il est enti\`{e}rement d\'{e}termin\'{e} par la
fonction analytique-r\'{e}elle $\mathcal{K}_{\Omega}$:%
\[
\mathcal{K}_{\Omega}(z)=K_{\Omega}(z,z)\qquad\left(  z\in\Omega\right)  ,
\]
qui est aussi appel\'{e}e noyau de Bergman de $\Omega$. De la m\^{e}me
mani\`{e}re, le noyau de Bergman \`{a} poids de $\left(  \Omega,p\right)  $
(noyau reproduisant de $H(\Omega,p)$) est not\'{e} $K_{\Omega,p}(z,t)$ et est
enti\`{e}rement d\'{e}termin\'{e} par $\mathcal{K}_{\Omega,p}(z)=K_{\Omega
,p}(z,z)$.

\begin{definition}
Soient $\Omega$ un domaine dans $V$ et $p:\Omega\rightarrow]0,+\infty\lbrack$
une fonction continue sur $\Omega$. On note $K_{\Omega,p^{k}}(z,w)$ (resp.
$\mathcal{K}_{\Omega,p^{k}}(z)$) le noyau de Bergman \`{a} poids de $\left(
\Omega,p^{k}\right)  $. On appelle \emph{noyau de Bergman virtuel} de $\left(
\Omega,p\right)  $ la fonction d\'{e}finie par%
\begin{equation}
L_{\Omega,p}\left(  z,w;r\right)  =L_{0}\left(  z,w;r\right)  =\sum
_{k=0}^{\infty}K_{\Omega,p^{k}}(z,w)r^{k}\quad(z,w\in\Omega,\ r\in\mathbb{C}).
\label{VB3a}%
\end{equation}
La fonction $\mathcal{L}_{0}(z,r)=L_{0}(z,z;r)$, d\'{e}finie par
\begin{equation}
\mathcal{L}_{\Omega,p}\left(  z,r\right)  =\mathcal{L}_{0}\left(  z,r\right)
=\sum_{k=0}^{\infty}\mathcal{K}_{\Omega,p^{k}}(z)r^{k}\quad(z\in\Omega
,\ r\geq0) \label{VB3}%
\end{equation}
est \'{e}galement appel\'{e}e noyau de Bergman virtuel de $\left(
\Omega,p\right)  $.
\end{definition}

\subsection{Noyau de Bergman de domaines de Hartogs}

Soit $\Omega\subset V$ et $p:\Omega\rightarrow]0,+\infty\lbrack$ une fonction
continue positive sur $\Omega$. On consid\`{e}re le domaine de Hartogs
$\widehat{\Omega}_{m}(p)$ au-dessus de $\Omega$, d\'{e}fini par
\[
\widehat{\Omega}_{m}(p)=\left\{  \left(  z,Z\right)  \in\Omega\times
\mathbb{C}^{m}\mid\left\Vert Z\right\Vert ^{2}<p(z)\right\}  .
\]
On munit ici $\mathbb{C}^{m}$ de la structure hermitienne standard et de la
forme volume associ\'{e}e
\[
\omega_{m}(Z)=\left(  \textstyle\frac{\operatorname*{i}}{2\pi}\partial
\overline{\partial}\left\Vert Z\right\Vert ^{2}\right)  ^{m}.
\]
Le domaine $\widehat{\Omega}_{m}(p)$ sera muni de la forme volume
\[
\omega_{V}(z)\wedge\omega_{m}(Z).
\]
Le th\'{e}or\`{e}me suivant montre comment calculer le noyau de Bergman des
domaines de Hartogs $\widehat{\Omega}_{m}(p)$ ($m>0$) \`{a} partir du noyau de
Bergman virtuel de $\left(  \Omega,p\right)  $.

\begin{theorem}
\label{TH1}Le noyau de Bergman $\widehat{K}_{m}$ (resp. $\widehat{\mathcal{K}%
}_{m}$) de $\widehat{\Omega}_{m}(p)$ est \'{e}gal \`{a}
\begin{align}
\widehat{K}_{m}\left(  (z,Z),(w,W)\right)   &  =L_{m}\left(  z,w;\left\langle
Z,W\right\rangle \right)  ,\label{VB5a}\\
\widehat{\mathcal{K}}_{m}(z,Z)  &  =\mathcal{L}_{m}\left(  z,\left\Vert
Z\right\Vert ^{2}\right)  ,\quad\label{VB5}%
\end{align}
o\`{u}
\begin{align}
L_{m}(z,w;r)  &  =\textstyle\frac{1}{m!}\frac{\partial^{m}}{\partial r^{m}%
}L_{0}(z,w;r),\label{VB12a}\\
\mathcal{L}_{m}(z,r)  &  =\textstyle\frac{1}{m!}\frac{\partial^{m}}{\partial
r^{m}}\mathcal{L}_{0}(z,r). \label{VB12}%
\end{align}

\end{theorem}

\subsection{Noyau de Bergman des domaines de Cartan-Hartogs}

Soit $\Omega$ un domaine homog\`{e}ne born\'{e} cercl\'{e} irr\'{e}ductible,
$N(z,t)$ sa norme g\'{e}n\'{e}rique, $\chi$ son polyn\^{o}me de Hua. Soit
\begin{equation}
\chi\left(  k\mu\right)  =\sum_{j=0}^{d}c_{j}\left(  \mu\right)
\textstyle\frac{\left(  k+1\right)  _{j}}{j!}=\sum_{j=0}^{d}\mu^{j}%
C_{d-j}\left(  \mu\right)  \left(  k+1\right)  _{j} \label{BergHart0}%
\end{equation}
la d\'{e}composition de $\chi\left(  k\mu\right)  $ suivant les factorielles
croissantes de $k$.

On consid\`{e}re le \emph{domaine de Cartan-Hartogs} $\widehat{\Omega}%
_{m}\left(  \mu\right)  $ construit au-dessus de $\Omega$ :
\begin{equation}
\widehat{\Omega}_{m}\left(  \mu\right)  =\left\{  \left(  z,Z\right)
\in\Omega\times\mathbb{C} ^{m},\quad\left\Vert Z\right\Vert ^{2}<N\left(
z,z\right)  ^{\mu}\right\}  . \label{hart1}%
\end{equation}
On note $\widehat{K}_{m,\mu}\left(  (z,Z),(w,W)\right)  $ le noyau de Bergman
de $\widehat{\Omega}_{m}\left(  \mu\right)  $ et
\[
\widehat{\mathcal{K}}_{m,\mu}(z,Z)=\widehat{K}_{m,\mu}\left(
(z,Z),(z,Z)\right)  .
\]
On note $L_{m,\mu}(z,w;r)$ le noyau de Bergman virtuel de $\left(
\Omega,N(z,z)^{\mu}\right)  $ et
\[
\mathcal{L}_{m,\mu}(z,r)=L_{m,\mu}(z,w;r).
\]

\begin{theorem}
Soit $\Omega$ un domaine homog\`{e}ne born\'{e} cercl\'{e} irr\'{e}ductible.
On a
\begin{align}
L_{0,\mu}\left(  z,w,r\right)   &  =\frac{C_{\Omega}}{N\left(  z,w\right)
^{g}}\sum_{j=0}^{d}c_{j}\left(  \mu\right)  \frac{1}{\left(  1-\xi\right)
^{j+1}},\label{BergHart1}\\
\mathcal{L}_{0,\mu}\left(  z,r\right)   &  =\frac{C_{\Omega}}{N\left(
z,z\right)  ^{g}}\sum_{j=0}^{d}c_{j}\left(  \mu\right)  \frac{1}{\left(
1-X\right)  ^{j+1}}, \label{BergHart2}%
\end{align}
o\`{u} $\xi$ et $X$ sont les fonctions d\'{e}finies par
\[
\xi(z,w,r)=\frac{r}{N(z,w)^{\mu}},\qquad X(z,r)=\frac{r}{N(z,z)^{\mu}}%
\]
et $C_{\Omega}=\frac{1}{\chi(0)\operatorname{vol}\Omega}$.
\end{theorem}

\begin{proof}
D'apr\`{e}s le lemme \ref{Bergman-poids} et (\ref{BergHart0}), on a
\begin{align*}
L_{0}\left(  z,w;r\right)   &  =\sum_{k=0}^{\infty}K_{\Omega,p^{k}}%
(z,w)r^{k}=C_{\Omega}\sum_{k=0}^{\infty}N(z,w)^{-g-k\mu}\chi(k\mu)r^{k}\\
&  =\frac{C_{\Omega}}{N(z,z)^{g}}\sum_{k=0}^{\infty}\chi(k\mu)\xi^{k}%
\end{align*}
o\`{u} $\xi=\frac{r}{N(z,w)^{\mu}}$. Si $P$ est un polyn\^{o}me
d\'{e}compos\'{e} sous la forme
\[
P(k)=\sum_{j=0}^{d}c_{j}\textstyle\frac{\left(  k+1\right)  _{j}}{j!},
\]
on a, pour $\left\vert \xi\right\vert <1$,%
\[
\sum_{k=0}^{\infty}P(k)\xi^{k}=\sum_{j=0}^{d}c_{j}\frac{1}{\left(
1-\xi\right)  ^{j+1}},
\]
d'o\`{u} le r\'{e}sultat pour $P(k)=\chi(k\mu)$ et $\chi\left(  k\mu\right)
=\sum_{j=0}^{d}c_{j}\left(  \mu\right)  \frac{\left(  k+1\right)  _{j}}{j!}$.
\end{proof}

\begin{notation}
On note $F_{\mu}=F_{\mu}^{0}$ la fonction rationnelle
\begin{equation}
F_{\mu}^{0}(\xi)=\sum_{j=0}^{d}c_{j}\left(  \mu\right)  \frac{1}{\left(
1-\xi\right)  ^{j+1}} \label{BergHart3}%
\end{equation}
et $P_{\mu}=P_{\mu}^{0}$ le polyn\^{o}me
\begin{equation}
P_{\mu}^{0}(\eta)=\sum_{j=0}^{d}c_{j}\left(  \mu\right)  \eta^{j}.
\label{BergHart4}%
\end{equation}
Plus g\'{e}n\'{e}ralement, pour $m$ entier positif, on note $F_{\mu}^{m}$ la
fraction rationnelle
\begin{equation}
F_{\mu}^{m}(\xi)=\sum_{j=0}^{d}\frac{(j+m)!}{j!}c_{j}\left(  \mu\right)
\frac{1}{\left(  1-\xi\right)  ^{j+m+1}} \label{BergHart5}%
\end{equation}
et $P_{\mu}^{m}$ le polyn\^{o}me%
\begin{equation}
P_{\mu}^{m}(\eta)=\frac{1}{m!}\sum_{j=0}^{d}\frac{(j+m)!}{j!}c_{j}\left(
\mu\right)  \eta^{j}. \label{BergHart6}%
\end{equation}

\end{notation}

\begin{theorem}
Le noyau de Bergman du domaine
\[
\widehat{\Omega}_{m}\left(  \mu\right)  =\left\{  \left(  z,Z\right)
\in\Omega\times\mathbb{C}^{m},\quad\left\Vert Z\right\Vert ^{2}<N\left(
z,z\right)  ^{\mu}\right\}
\]
est%
\begin{equation}
\widehat{K}_{m,\mu}\left(  (z,Z),(w,W)\right)  =\frac{1}{m!}\frac{C}{N\left(
z,w\right)  ^{g+m\mu}}F_{\mu}^{m}(\xi), \label{BergHart7}%
\end{equation}
o\`{u} $\xi:\widehat{\Omega}_{m}\left(  \mu\right)  \times\widehat{\Omega}%
_{m}\left(  \mu\right)  \rightarrow\mathbb{C}$ est d\'{e}finie par%
\begin{equation}
\xi\left(  (z,Z),(w,W)\right)  =\frac{\left\langle Z,W\right\rangle
}{N(z,w)^{\mu}}. \label{BergHart8}%
\end{equation}

\end{theorem}

\begin{proof}
On a
\begin{align*}
L_{0,\mu}\left(  z,w,r\right)   &  =\frac{C}{N\left(  z,w\right)  ^{g}}%
\sum_{j=0}^{d}c_{j}\left(  \mu\right)  \frac{1}{\left(  1-\xi\right)  ^{j+1}%
},\\
\xi(z,w,r)  &  =\frac{r}{N(z,w)^{\mu}}%
\end{align*}
et
\begin{align*}
L_{m,\mu}\left(  z,w,r\right)   &  =\frac{1}{m!}\frac{\partial^{m}}{\partial
r^{m}}L_{0,\mu}\left(  z,w,r\right) \\
&  =\frac{1}{m!}\frac{C}{N\left(  z,w\right)  ^{g+m\mu}}\sum_{j=0}^{d}%
\frac{(j+m)!}{j!}c_{j}\left(  \mu\right)  \frac{1}{\left(  1-\xi\right)
^{j+m+1}}.
\end{align*}
On en d\'{e}duit le noyau de Bergman de $\widehat{\Omega}_{m}\left(
\mu\right)  $ :
\[
\widehat{K}_{m,\mu}\left(  (z,Z),(w,W)\right)  =L_{m,\mu}\left(
z,w;\left\langle Z,W\right\rangle \right)  =\frac{1}{m!}\frac{C}{N\left(
z,w\right)  ^{g+m\mu}}F_{\mu}^{m}(\xi).
\]
\end{proof}

Soit $\eta:\widehat{\Omega}_{m}\left(  \mu\right)  \times\widehat{\Omega}%
_{m}\left(  \mu\right)  \rightarrow\mathbb{C}$ d\'{e}finie par
\begin{equation}
\eta\left(  (z,Z),(w,W)\right)  =\frac{1}{1-\xi\left(  (z,Z),(w,W)\right)
}=\frac{N(z,w)^{\mu}-\left\langle Z,W\right\rangle }{N(z,w)^{\mu}}.
\label{BergHart9}%
\end{equation}
Le noyau de Bergman de $\widehat{\Omega}_{m}\left(  \mu\right)  $ s'\'{e}crit
encore
\begin{equation}
\widehat{K}_{m,\mu}\left(  (z,Z),(w,W)\right)  =\frac{C}{N\left(  z,w\right)
^{g+m\mu}}\eta^{m+1}P_{\mu}^{m}(\eta). \label{BergHart7a}%
\end{equation}

\begin{definition}
Le polyn\^{o}me
\[
P_{\mu}^{m}(\eta)=\frac{1}{m!}\sum_{j=0}^{d}\frac{(j+m)!}{j!}c_{j}\left(
\mu\right)  \eta^{j}=\sum_{j=0}^{d}(m+1)_{j}C_{d-j}(\mu)\mu^{j}\eta^{j}%
\]
sera appel\'{e} \emph{polyn\^{o}me repr\'{e}sentatif du noyau de Bergman de}
$\widehat{\Omega}_{m}\left(  \mu\right)  $.
\end{definition}

\section{\label{LuProblem}Probl\`{e}me de Lu Qikeng pour les domaines de
Cartan-Hartogs}

Le \emph{probl\`{e}me de Lu Qikeng} pour un ouvert $U$ de $\mathbb{C}^{n}$
consiste \`{a} d\'{e}terminer si le noyau de Bergman $K_{U}(z,w)$ de ce
domaine peut avoir des z\'{e}ros. Un domaine $U$ est appel\'{e} \emph{domaine
de Lu Qikeng} si son noyau de Bergman ne s'annule pas dans $U\times U$.

Soit $\Omega$ un domaine homog\`{e}ne born\'{e} cercl\'{e} irr\'{e}ductible,
$N(z,t)$ sa norme g\'{e}\-n\'{e}\-rique, $\chi$ son polyn\^{o}me de Hua. Pour
$\mu>0$ et $m$ entier positif, on consid\`{e}re le \emph{domaine de
Cartan-Hartogs} $\widehat{\Omega}_{m}\left(  \mu\right)  $ construit au-dessus
de $\Omega$ :
\[
\widehat{\Omega}_{m}\left(  \mu\right)  =\left\{  \left(  z,Z\right)
\in\Omega\times\mathbb{C}^{m},\quad\left\Vert Z\right\Vert ^{2}<N\left(
z,z\right)  ^{\mu}\right\}  .
\]
Le noyau de Bergman de $\widehat{\Omega}_{m}\left(  \mu\right)  $ s'\'{e}crit
\begin{equation}
\widehat{K}_{m,\mu}\left(  (z,Z),(w,W)\right)  =\frac{C}{N\left(  z,w\right)
^{g+m\mu}}\eta^{m+1}P_{\mu}^{m}(\eta), \label{Lu3}%
\end{equation}
o\`{u} le polyn\^{o}me $P_{\mu}^{m}$ est le \emph{polyn\^{o}me
repr\'{e}sentatif du noyau de Bergman}, d\'{e}fini \`{a} partir du
polyn\^{o}me de Hua $\chi$ par les relations
\begin{align}
\chi\left(  k\mu\right)   &  =\sum_{j=0}^{d}\mu^{j}C_{d-j}\left(  \mu\right)
\left(  k+1\right)  _{j},\label{Lu4}\\
P_{\mu}^{m}(\eta)  &  =\sum_{j=0}^{d}(m+1)_{j}\mu^{j}C_{d-j}(\mu)\eta^{j},
\label{Lu5}%
\end{align}
et la fonction $\eta$ est d\'{e}finie par
\begin{align}
\xi\left(  (z,Z),(w,W)\right)   &  =\frac{\left\langle Z,W\right\rangle
}{N(z,w)^{\mu}},\label{Lu1}\\
\eta\left(  (z,Z),(w,W)\right)   &  =\frac{1}{1-\xi\left(  (z,Z),(w,W)\right)
}. \label{Lu2}%
\end{align}

\begin{lemma}
Soient $\xi$ et $\eta$ les fonctions d\'{e}finies sur $\widehat{\Omega}%
_{m}\left(  \mu\right)  \times\widehat{\Omega}_{m}\left(  \mu\right)  $ par
(\ref{Lu1}) et (\ref{Lu2}). Alors l'image de $\xi$ est le disque unit\'{e}
$\Delta$ de $\mathbb{C}$ et l'image de $\eta$ est le demi-plan $\left\{
\operatorname{Re}\eta>\frac{1}{2}\right\}  $.
\end{lemma}

\begin{proof}
Le noyau de Bergman de $\Omega$ est
\[
K(z,t)=C_{0}N(z,t)^{-g},
\]
avec $C_{0}=(\operatorname{vol}\Omega)^{-1}$. De la propri\'{e}t\'{e} connue
du noyau de Bergman :
\begin{equation}
\left\vert K(z,t)\right\vert ^{2}\leq K(z,z)K(t,t), \label{inegal2}%
\end{equation}
on d\'{e}duit
\begin{equation}
\left\vert N(z,t)\right\vert ^{2}\geq N(z,z)N(t,t)\qquad(z,t\in\Omega).
\label{inegal1}%
\end{equation}
Dans $\widehat{\Omega}_{m}(\mu)\times\widehat{\Omega}_{m}(\mu)$, on a donc
\[
\left\vert \xi\left(  (z,Z),(t,T)\right)  \right\vert ^{2}=\frac{\left\vert
\left\langle Z,T\right\rangle \right\vert ^{2}}{N(z,t)^{2\mu}}\leq
\frac{\left\Vert Z\right\Vert ^{2}}{N(z,z)^{\mu}}\frac{\left\Vert T\right\Vert
^{2}}{N(t,t)^{\mu}}<1.
\]
La fonction $\xi$ prend donc ses valeurs dans le disque unit\'{e} $\Delta$ de
$\mathbb{C}$. Comme%
\[
\xi\left(  (0,Z),(0,\operatorname{e}^{\operatorname{i}\theta}Z)\right)
=\operatorname{e}^{\operatorname{i}\theta}\left\Vert Z\right\Vert ^{2},
\]
l'image de $\xi$ est \'{e}gale \`{a} $\Delta$.
\end{proof}

L'expression (\ref{Lu3}) du noyau de Bergman de $\widehat{\Omega}_{m}\left(
\mu\right)  $ entra\^{\i}ne alors imm\'{e}\-dia\-te\-ment

\begin{theorem}
\label{SolAlgPbLu}Le domaine $\widehat{\Omega}_{m}\left(  \mu\right)  $ est un
domaine de Lu Qikeng si et seulement si le polyn\^{o}me $P_{\mu}^{m}$ ne
s'annule pas dans $\left\{  \operatorname{Re}\eta>\frac{1}{2}\right\}  $.
\end{theorem}

Le probl\`{e}me de Lu Qikeng pour les domaines $\widehat{\Omega}_{m}\left(
\mu\right)  $ est ainsi ramen\'{e} \`{a} la localisation des racines du
polyn\^{o}me $P_{\mu}^{m}$, qui a le m\^{e}me degr\'{e} que le polyn\^{o}me de
Hua $\chi$ et s'en d\'{e}duit alg\'{e}briquement.

\section{\label{LowDim}Solution du probl\`{e}me de Lu Qikeng pour une base de
faible dimension}

Dans cette section, nous donnons la solution compl\`{e}te du probl\`{e}me de
Lu Qikeng pour les domaines de Cartan-Hartogs $\widehat{\Omega}_{m}\left(
\mu\right)  $, lorsque la base $\Omega$ est un domaine born\'{e}
sym\'{e}trique irr\'{e}ductible de dimension au plus $4$.

\subsection{Cas o\`{u} $\Omega$ est le disque unit\'{e} de $\mathbb{C}$}

Le domaine $\Omega$ est le disque unit\'{e} de $\mathbb{C}$ et le domaine
$\widehat{\Omega}_{m}\left(  \mu\right)  $ est
\begin{align*}
\widehat{\Omega}_{m}\left(  \mu\right)   &  =\left\{  \left(  z,Z\right)
\in\Omega\times\mathbb{C}^{m},\quad\left\Vert Z\right\Vert ^{2}<\left(
1-\left\vert z\right\vert ^{2}\right)  ^{\mu}\right\} \\
&  =\left\{  \left(  z,Z\right)  \in\mathbb{C}\times\mathbb{C}^{m}%
,\quad\left\vert z\right\vert ^{2}+\left\Vert Z\right\Vert ^{2/\mu}<1\right\}
.
\end{align*}
On a (voir l'annexe \ref{TableI-1-1})%
\[
\textstyle P_{\mu}^{m}\left(  \frac{1}{2}\right)  =\frac{(m-1)\mu}{2}+1,
\]
qui est positif pour tout $\mu>0$ et tout $m\geq1$. La racine de $P_{\mu}^{m}$
est donc toujours inf\'{e}rieure \`{a} $\frac{1}{2}$ ; d'o\`{u}

\begin{theorem}
\label{SolPbLuI-1-1}Si $\Omega$ est le disque unit\'{e} de $\mathbb{C}$, le
domaine $\widehat{\Omega}_{m}\left(  \mu\right)  $ est un domaine de Lu Qikeng
pour tout $\mu>0$ et tout entier $m\geq1$.
\end{theorem}

Ce r\'{e}sultat facile et connu (voir par exemple \cite{Yin2006}) n'est
cit\'{e} ici que pour comparaison avec les r\'{e}sultats qui suivront et parce
qu'il illustre la m\'{e}thode employ\'{e}e lorsque $\Omega$ est de dimension
$>1$.

\subsection{Cas o\`{u} $\Omega$ est de type $I_{1,2}$}

Le domaine $\Omega$ est la boule unit\'{e} de $\mathbb{C}^{2}$ et le domaine
$\widehat{\Omega}_{m}\left(  \mu\right)  $ est
\begin{align*}
\widehat{\Omega}_{m}\left(  \mu\right)   &  =\left\{  \left(  z,Z\right)
\in\Omega\times\mathbb{C}^{m},\quad\left\Vert Z\right\Vert ^{2}<\left(
1-\left\Vert z\right\Vert ^{2}\right)  ^{\mu}\right\} \\
&  =\left\{  \left(  z,Z\right)  \in\mathbb{C}^{2}\times\mathbb{C}^{m}%
,\quad\left\Vert z\right\Vert ^{2}+\left\Vert Z\right\Vert ^{2/\mu}<1\right\}
.
\end{align*}
Le polyn\^{o}me de Hua est
\[
\chi\left(  s\right)  =\left(  s+1\right)  \left(  s+2\right)  .
\]
Le polyn\^{o}me repr\'{e}sentatif du noyau de Bergman de $\widehat{\Omega}%
_{m}\left(  \mu\right)  $ est (voir annexe \ref{TableI-1-2})%
\[
P_{\mu}^{m}(\eta)=(1-\mu)(2-\mu)+3(m+1)\mu(1-\mu)\eta+(m+1)(m+2)\mu^{2}%
\eta^{2}.
\]
D'apr\`{e}s le th\'{e}or\`{e}me \ref{SolAlgPbLu}, le domaine $\widehat{\Omega
}_{m}\left(  \mu\right)  $ est un domaine de Lu Qikeng si et seulement si le
polyn\^{o}me $P_{\mu}^{m}(\eta)$ a toutes ses racines dans le demi plan
$\left\{  \operatorname{Re}\left(  z\right)  \leq\frac{1}{2}\right\}  $.
D'apr\`{e}s la proposition \ref{LocDegre2} ces racines sont dans le demi plan
$\left\{  \operatorname{Re}\left(  z\right)  <\frac{1}{2}\right\}  $ si et
seulement si
\begin{align}
P_{\mu}^{m}\left(  \textstyle\frac{1}{2}\right)   &  >0,\label{A1}\\
\textstyle\frac{\operatorname*{d}P_{\mu}^{m}}{\operatorname*{d}\eta}\left(
\frac{1}{2}\right)   &  >0. \label{A2}%
\end{align}
On a (annexe \ref{TableI-1-2})%
\[
\textstyle\frac{1}{(m+1)\mu}\frac{\operatorname*{d}P_{\mu}^{m}}%
{\operatorname*{d}\eta}\left(  \frac{1}{2}\right)  =3+(m-1)\mu;
\]
la condition (\ref{A2}) est donc v\'{e}rifi\'{e}e pour tout $\mu>0$ et tout
$m\geq1$.

On a d'autre part (annexe \ref{TableI-1-2})%
\[
\textstyle P_{\mu}^{m}\left(  \frac{1}{2}\right)  =\frac{m(m-3)}{4}\mu
^{2}+\frac{3(m-1)}{2}\mu+2.
\]
Si $m\geq3$, tous les coefficients de ce polyn\^{o}me en $\mu$ sont non
n\'{e}gatifs et on a $P_{\mu}^{m}\left(  \frac{1}{2}\right)  >0$ pour tout
$\mu>0$.

Si $m=1$, on a
\[
\textstyle P_{\mu}^{1}\left(  \frac{1}{2}\right)  =2-\frac{\mu^{2}}{2},
\]
qui est strictement positif si et seulement si $\mu<\mu_{1}=2$.

Si $m=2$, on a
\[
\textstyle P_{\mu}^{2}\left(  \frac{1}{2}\right)  =2+\frac{3}{2}\mu-\frac
{1}{2}\mu^{2},
\]
qui est strictement positif si et seulement si $\mu<\mu_{2}=4$. Les conditions
(\ref{A2}) et (\ref{A1}) sont donc v\'{e}rifi\'{e}es pour $\mu<\mu_{m}$,
$\mu_{m}$ \'{e}tant la racine positive de $P_{\mu}^{m}\left(  \frac{1}%
{2}\right)  $. Comme l'ensemble des racines varie contin\^{u}ment en fonction
de $\mu$, le domaine est de Lu Qikeng si et seulement si $\mu\leq\mu_{m}$. En
conclusion :

\begin{theorem}
\label{SolPbLuI-1-2}Soit $\Omega$ la boule unit\'{e} de $\mathbb{C}^{2}$. Le
domaine $\widehat{\Omega}_{m}\left(  \mu\right)  $ est un domaine de Lu Qikeng

\begin{itemize}
\item pour $m=1$, si et seulement si $\mu\leq2$ ;

\item pour $m=2$, si et seulement si $\mu\leq4$ ;

\item pour $m\geq3$, quel que soit $\mu>0$.
\end{itemize}
\end{theorem}

Si $m=1$, ce r\'{e}sultat est d\^{u} \`{a} H.P. Boas, Siqi Fu, E. Straube
(\cite{BoasFuStraube1999}) ; voir aussi \cite{Yin2006}. Pour $m>1$, les
r\'{e}sultats du th\'{e}or\`{e}me sont nouveaux.

\subsection{Cas o\`{u} $\Omega$ est de type $I_{1,3}$}

Le domaine $\Omega$ est la boule hermitienne de dimension $3$. Le polyn\^{o}me
de Hua est
\[
\chi\left(  s\right)  =\left(  s+1\right)  \left(  s+2\right)  (s+3).
\]
Le polyn\^{o}me repr\'{e}sentatif du noyau de Bergman de $\widehat{\Omega}%
_{m}\left(  \mu\right)  $ est (voir annexe \ref{TableI-1-3})
\begin{align*}
P_{\mu}^{m}(\eta)   =(&1-\mu)(2-\mu)\left(  3-\mu\right)  +(m+1)\mu
(1-\mu)\left(  11-7\mu\right)  \eta\\
&  +6(m+1)(m+2)\left(  1-\mu\right)  \mu^{2}\eta^{2}+\left(  m+1\right)
\left(  m+2\right)  \left(  m+3\right)  \mu^{3}\eta^{3}.
\end{align*}
On a
\[
\textstyle P_{\mu}^{m}\left(  \frac{1}{2}\right)  =6+\frac{11\left(
m-1\right)  }{2}\mu+\frac{3m(m-3)}{2}\mu^{2}+\frac{\left(  m-1\right)  \left(
m^{2}-5m-2\right)  }{8}\mu^{3}.
\]
Pour $1\leq m\leq5$, soit $\mu_{m}$ l'unique racine positive du polyn\^{o}me
$q_{m}$ d\'{e}fini par $q_{m}(\mu)=P_{\mu}^{m}\left(  \frac{1}{2}\right)  $.
On a (cf. proposition \ref{Signe-I-1-3(0)})%
\[
0<\mu_{1}=\sqrt{2}<\mu_{2}<\mu_{3}<\mu_{4}<\mu_{5}.
\]

Du th\'{e}or\`{e}me \ref{Racines-I-1-3}, on d\'{e}duit la solution
compl\`{e}te du probl\`{e}me de Lu Qikeng lorsque $\Omega$ est de type
$I_{1,3}$ :

\begin{theorem}
\label{SolPbLuI-1-3}Soit $\Omega$ une boule hermitienne de dimension $3$. Si
$m\geq6$, le domaine de Cartan-Hartogs $\widehat{\Omega}_{m}\left(
\mu\right)  $ est un domaine de Lu Qikeng pour tout $\mu\in]0,+\infty\lbrack$.
Si $1\leq m\leq5$, le domaine de Cartan-Hartogs $\widehat{\Omega}_{m}\left(
\mu\right)  $ est un domaine de Lu Qikeng si et seulement si $0<\mu\leq\mu
_{m}.$
\end{theorem}

Pour $m=1$, ce r\'{e}sultat a \'{e}t\'{e} obtenu par Weiping Yin
(\cite{Yin2006}) par une m\'{e}thode diff\'{e}rente mais essentiellement
\'{e}quivalente. Les r\'{e}sultats du th\'{e}or\`{e}me sont nouveaux pour
$m>1$.

\subsection{Cas o\`{u} $\Omega$ est de type $IV_{3}$}

Le domaine $\Omega$ est la boule de Lie de dimension $3$ (isomorphe au domaine
sym\'{e}trique associ\'{e} \`{a} l'espace $\mathcal{S}_{2}(\mathbb{C})$ des
matrices sym\'{e}triques $(2,2)$). Les invariants num\'{e}riques sont $a=1$,
$b=0$, $r=2$. Le polyn\^{o}me de Hua est
\[
\textstyle\chi\left(  s\right)  =\left(  s+1\right)  \left(  s+\frac{3}%
{2}\right)  \left(  s+2\right)  .
\]
Le polyn\^{o}me repr\'{e}sentatif du noyau de Bergman de $\widehat{\Omega}%
_{m}\left(  \mu\right)  $ est (voir annexe \ref{TableIV-3})%
\begin{align*}
P_{\mu}^{m}   (\eta) =\textstyle(&1-\mu)(2-\mu)\left(  \frac{3}{2}-\mu\right)
+(m+1)\mu(1-\mu)\left(  \frac{13}{2}-7\mu\right)  \eta\\
&  \quad+\textstyle3(m+1)(m+2)\left(  \frac{3}{2}-2\mu\right)  \mu^{2}\eta
^{2}+\left(  m+1\right)  \left(  m+2\right)  \left(  m+3\right)  \mu^{3}%
\eta^{3}.
\end{align*}
On a
\[
\textstyle P_{\mu}^{m}\left(  \frac{1}{2}\right)  =\frac{\left(  m-1\right)
\left(  m^{2}-5m-2\right)  }{8}\mu^{3}+\frac{9m\left(  m-3\right)  }{8}\mu
^{2}+\frac{13\left(  m-1\right)  }{4}\mu+3.
\]
Pour $1\leq m\leq5$, soient $q_{m}$ les polyn\^{o}mes d\'{e}finis par
\[
\textstyle q_{m}(\mu)=P_{\mu}^{m}\left(  \frac{1}{2}\right)
\]
et soit $\mu_{m}$ l'unique racine positive du polyn\^{o}me $q_{m}$. On a (cf.
proposition \ref{Signe-IV-3(0)})%
\[
\textstyle0<\mu_{1}=\frac{2}{\sqrt{3}}<\mu_{2}<\mu_{3}<\mu_{4}<\mu_{5}.
\]

Du th\'{e}or\`{e}me \ref{Racines-IV-3}, on d\'{e}duit la solution compl\`{e}te
du probl\`{e}me de Lu Qikeng lorsque $\Omega$ est de type $IV_{3}$ :

\begin{theorem}
\label{SolPbLuIV-3}Soit $\Omega$ une boule de Lie de dimension $3$. Si
$m\geq6$, le domaine de Cartan-Hartogs $\widehat{\Omega}_{m}\left(
\mu\right)  $ est un domaine de Lu Qikeng pour tout $\mu\in]0,+\infty\lbrack$.
Si $1\leq m\leq5$, le domaine de Cartan-Hartogs $\widehat{\Omega}_{m}\left(
\mu\right)  $ est un domaine de Lu Qikeng si et seulement si $0<\mu\leq\mu
_{m}.$
\end{theorem}

Les r\'{e}sultats de ce th\'{e}or\`{e}me sont enti\`{e}rement nouveaux.

\subsection{Cas o\`{u} $\Omega$ est de type $I_{1,4}$}

Le domaine $\Omega$ est la boule hermitienne de dimension $4$. Le polyn\^{o}me
de Hua est
\[
\chi\left(  s\right)  =\left(  s+1\right)  \left(  s+2\right)  (s+3)\left(
s+4\right)  .
\]
Le polyn\^{o}me repr\'{e}sentatif du noyau de Bergman de $\widehat{\Omega}%
_{m}\left(  \mu\right)  $ est (voir annexe \ref{TableI-1-4})
\begin{align*}
P_{\mu}^{m}(\eta)=(  &  1-\mu)(2-\mu)\left(  3-\mu\right)  (4-\mu)\\
&  +5(m+1)(1-\mu)(5-3\mu)(2-\mu)\mu\eta\\
&  +5(m+1)_{2}\left(  1-\mu\right)  (7-5\mu)\mu^{2}\eta^{2}\\
&  +10\left(  m+1\right)  _{3}(1-\mu)\mu^{3}\eta^{3}+\left(  m+1\right)
_{4}\mu^{4}\eta^{4}%
\end{align*}
On a
\begin{align*}
\textstyle P_{\mu}^{m}\left(  \frac{1}{2}\right)  =(4+m\mu)  &  \left[
\textstyle6+\frac{1}{4}(19m-25)\mu+m(m-5)\mu^{2}\right. \\
&  \quad+\left.  \textstyle\frac{1}{16}(m^{3}-10m^{2}+15m+10)\mu^{3}\right]  .
\end{align*}
Pour $1\leq m\leq7$, soit $\mu_{m}$ la plus petite racine positive du
polyn\^{o}me $q_{m}$ d\'{e}fini par $q_{m}(\mu)=P_{\mu}^{m}\left(  \frac{1}%
{2}\right)  $ . On a (cf. proposition \ref{Signe-I-1-4(0)})%
\[
\textstyle0<\mu_{1}=\sqrt{\frac{3}{2}}<\mu_{2}<\mu_{3}<\mu_{4}<\mu_{5}<\mu
_{6}<\mu_{7}.
\]

Du th\'{e}or\`{e}me (\ref{LocTypeI-1-4}), on d\'{e}duit la solution
compl\`{e}te du probl\`{e}me de Lu Qikeng lorsque $\Omega$ est de type
$I_{1,4}$ :

\begin{theorem}
\label{SolPbLuI-1-4}Soit $\Omega$ une boule hermitienne de dimension $4$. Si
$m\geq8$, le domaine de Cartan-Hartogs $\widehat{\Omega}_{m}\left(
\mu\right)  $est un domaine de Lu Qikeng pour tout $\mu\in]0,+\infty\lbrack$.
Si $1\leq m\leq7$, le domaine $\widehat{\Omega}_{m}\left(  \mu\right)  $ est
un domaine de Lu Qikeng si et seulement si $0<\mu\leq\mu_{m}.$
\end{theorem}

Pour $m=1$, ce r\'{e}sultat a \'{e}t\'{e} obtenu par Liyou Zhang et Jong-do
Park (2006, non publi\'{e}). Les r\'{e}sultats du th\'{e}or\`{e}me sont
nouveaux pour $m>1$.

\subsection{Cas o\`{u} $\Omega$ est de type $IV_{4}$}

Le domaine $\Omega$ est la boule de Lie de dimension $4$. Le polyn\^{o}me de
Hua est
\[
\chi\left(  s\right)  =\left(  s+1\right)  \left(  s+2\right)  ^{2}(s+3).
\]
Le polyn\^{o}me repr\'{e}sentatif du noyau de Bergman de $\widehat{\Omega}%
_{m}\left(  \mu\right)  $ est (voir annexe \ref{TableIV-4})%
\begin{align*}
P_{\mu}^{m}(\eta) =(  &  1-\mu)\left(  2-\mu\right)  ^{2}(3-\mu)+(m+1)(1-\mu
)(7-5\mu)(4-3\mu)\mu\eta\\
&  \ +(m+1)_{2}(1-\mu)(23-25\mu)\mu^{2}\eta^{2}+2\left(  m+1\right)
_{3}(4-5\mu)\mu^{3}\eta^{3}\\
&  \ +\left(  m+1\right)  _{4}\mu^{4}\eta^{4}.
\end{align*}
On a%
\begin{align*}
\textstyle P_{\mu}^{m}\left(  \frac{1}{2}\right)     =1&2+14(m-1)\mu
+\textstyle\frac{23m(m-3)}{4}\mu^{2}\\
&  +(m-1)(m^{2}-5m-2)\mu^{3}+\textstyle\frac{m^{3}-10m^{2}+15m+10}{16}\mu^{4}.
\end{align*}

Pour $1\leq m\leq7$, soit $\mu_{m}$ la plus petite racine positive du
polyn\^{o}me $q_{m}$ d\'{e}fini par $q_{m}(\mu)=P_{\mu}^{m}\left(  \frac{1}%
{2}\right)  $ On a (cf. proposition \ref{Signe-IV-4(0)})%
\[
\textstyle0<\mu_{1}=\frac{1}{2}\sqrt{23-\sqrt{337}}<\mu_{2}<\mu_{3}<\mu
_{4}<\mu_{5}<\mu_{6}<\mu_{7}.
\]

Du th\'{e}or\`{e}me (\ref{LocTypeIV-4}), on d\'{e}duit la solution
compl\`{e}te du probl\`{e}me de Lu Qikeng lorsque $\Omega$ est de type
$IV_{4}$ :

\begin{theorem}
\label{SolPbLuIV-4}Soit $\Omega$ une boule de Lie de dimension $4$. Si
$m\geq8$, le domaine de Cartan-Hartogs $\widehat{\Omega}_{m}\left(
\mu\right)  $ est un domaine de Lu Qikeng pour tout $\mu\in]0,+\infty\lbrack$.
Si $1\leq m\leq7$, le domaine $\widehat{\Omega}_{m}\left(  \mu\right)  $ est
un domaine de Lu Qikeng si et seulement si $0<\mu\leq\mu_{m}$.
\end{theorem}

Les r\'{e}sultats de ce th\'{e}or\`{e}me sont enti\`{e}rement nouveaux.

\appendix

\section{\label{Localisation}Localisation des racines}

Le probl\`{e}me de Lu Qikeng pour les domaines de Cartan-Hartogs
$\widehat{\Omega}_{m}\left(  \mu\right)  $ a \'{e}t\'{e} ramen\'{e} \`{a} la
localisation des racines du polyn\^{o}me $P_{\mu}^{m}$ par rapport au
demi-plan $\left\{  \operatorname{Re}z>\frac{1}{2}\right\}  $
(th\'{e}or\`{e}me \ref{SolAlgPbLu}). Nous rappelons ci-dessous le
\emph{crit\`{e}re de Routh--Hurwitz} qui permet de d\'{e}terminer quand un
polyn\^{o}me \`{a} coefficients r\'{e}els a toutes ses racines dans le
demi-plan $\left\{  \operatorname{Re}z<0\right\}  $. Dans les paragraphes
suivants, nous en d\'{e}duisons les crit\`{e}res utilis\'{e}s pour
r\'{e}soudre le probl\`{e}me de Lu Qikeng pour un domaine de Cartan--Hartogs
dont la base $\Omega$ est de dimension au plus $4$.

\subsection{Crit\`{e}re de Routh--Hurwitz}

Soit $P$ un polyn\^{o}me \`{a} coefficients r\'{e}els de degr\'{e} $n$
\begin{equation}
P\left(  z\right)  =a_{0}z^{n}+a_{1}z^{n-1}+\ldots+a_{n-1}z+a_{n};
\label{CRH1}%
\end{equation}
on suppose $a_{0}>0$. Le polyn\^{o}me $P$ est dit \emph{stable} si toutes ses
racines ont des parties r\'{e}elles n\'{e}gatives.

L'\'{e}tude de la stabilit\'{e} des polyn\^{o}mes intervient dans la
th\'{e}orie du contr\^{o}le. En effet, l'\'{e}quation caract\'{e}ristique d'un
syst\`{e}me d'\'{e}quations diff\'{e}rentielles li\-n\'{e}ai\-res est un
polyn\^{o}me et la stabilit\'{e} du syst\`{e}me se traduit par le fait que
toutes les racines de l'\'{e}quation caract\'{e}ristique soient dans le
demi-plan n\'{e}gatif.

En 1875, le m\'{e}canicien anglais Routh a \'{e}labor\'{e} un algorithme qui
permet de d\'{e}terminer si un polyn\^{o}me est stable (et plus
g\'{e}n\'{e}ralement de localiser ses racines par rapport \`{a} $\left\{
\operatorname{Re}z=0\right\}  $). Vingt ans plus tard, le math\'{e}maticien
allemand Hurwitz a donn\'{e} un crit\`{e}re \'{e}quivalent, dans une forme
diff\'{e}rente, utilisant les d\'{e}terminants (\emph{d\'{e}terminants de
Hurwitz})
\begin{align*}
\Delta_{1}  &  =a_{1},\quad\Delta_{2}=%
\begin{vmatrix}
a_{1} & a_{3}\\
a_{0} & a_{2}%
\end{vmatrix}
,\quad\Delta_{3}=
\begin{vmatrix}
a_{1} & a_{3} & a_{5}\\
a_{0} & a_{2} & a_{4}\\
0 & a_{1} & a_{3}%
\end{vmatrix}
,\\
\Delta_{n}  &  =
\begin{vmatrix}
a_{1} & a_{3} & a_{5} & \cdots & 0\\
a_{0} & a_{2} & a_{4} & \cdots & 0\\
0 & a_{1} & a_{3} & \cdots & 0\\
\vdots & a_{0} & a_{2} & \cdots & \cdots\\
0 & 0 & 0 & \cdots & a_{n}%
\end{vmatrix}
.
\end{align*}

\begin{criterion}
[Crit\`{e}re de Routh--Hurwitz]Toutes les racines du polyn\^{o}me (\ref{CRH1})
ont des parties r\'{e}elles n\'{e}gatives si et seulement si les $n$
d\'{e}terminants de Hurwitz sont positifs :
\[
\Delta_{1}>0,\quad\Delta_{2}>0,\quad\cdots,\quad\Delta_{n}>0.
\]

\end{criterion}

En 1914, les math\'{e}maticiens fran\c{c}ais Li\'{e}nard et Chipart ont
\'{e}tabli un crit\`{e}re de stabilit\'{e}, diff\'{e}rent de celui de
Routh-Hurwitz, en constatant que lorsque les coefficients $a_{0}$, $a_{1}$ ,
$\cdots$ , $a_{n}$ sont positifs, les conditions sur la positivit\'{e} des
$\Delta_{i}$ ne sont pas ind\'{e}pendantes. Par exemple, pour $n=4$, les
conditions de stabilit\'{e} se r\'{e}duisent \`{a} $a_{1}>0$, $a_{2}>0$,
$a_{4}>0$ et $\Delta_{3}>0$.

\subsection{Polyn\^{o}mes de degr\'{e} $2$}

Soit
\[
Q\left(  z\right)  =a_{0}z^{2}+a_{1}z+a_{2},
\]
avec $a_{0}>0$. Les d\'{e}terminants de Hurwitz sont
\[
\Delta_{1}=a_{1},\quad\Delta_{2}=%
\begin{vmatrix}
a_{1} & 0\\
a_{0} & a_{2}%
\end{vmatrix}
=a_{1}a_{2}.
\]
Le polyn\^{o}me $Q$ est stable si et seulement si on a
\begin{equation}
a_{1}>0,\quad a_{2}>0. \label{RouthDeg2}%
\end{equation}

\begin{proposition}
\label{LocDegre2}Le polyn\^{o}me du second degr\'{e} $P$ a toutes ses racines
dans $\left\{  \operatorname{Re}z<\frac{1}{2}\right\}  $ si et seulement si le
polyn\^{o}me $Q(z)=P\left(  \frac{1}{2}+z\right)  $ est stable,
c'est-\`{a}-dire si
\begin{equation}
\textstyle P\left(  \frac{1}{2}\right)  >0,\quad P^{\prime}\left(  \frac{1}%
{2}\right)  >0. \label{CritDeg2}%
\end{equation}

\end{proposition}

Il est \'{e}galement facile d'\'{e}tablir cette proposition directement, sans
recours au crit\`{e}re de Routh--Hurwitz.

\subsection{Polyn\^{o}mes de degr\'{e} $3$}

Soit
\[
Q\left(  z\right)  =a_{0}z^{3}+a_{1}z^{2}+a_{2}z+a_{3},
\]
avec $a_{0}>0$. Les d\'{e}terminants de Hurwitz sont
\[
\Delta_{1}=a_{1},\quad\Delta_{2}=%
\begin{vmatrix}
a_{1} & a_{3}\\
a_{0} & a_{2}%
\end{vmatrix}
,\quad\Delta_{3}=%
\begin{vmatrix}
a_{1} & a_{3} & 0\\
a_{0} & a_{2} & 0\\
0 & a_{1} & a_{3}%
\end{vmatrix}
=a_{3}\Delta_{2}.
\]
Le polyn\^{o}me $Q$ est stable si et seulement si on a
\[
a_{1}>0,\quad\Delta_{2}>0,\quad a_{3}>0.
\]
Comme $\Delta_{2}=a_{1}a_{2}-a_{0}a_{3}$, cette condition \'{e}quivaut \`{a}
\begin{equation}
a_{3}>0,\quad a_{2}>0,\quad\Delta_{2}>0. \label{LCDeg3}%
\end{equation}

Soit $P(z)=\alpha+\beta z+\gamma z^{2}+\delta z^{3}$ un polyn\^{o}me de
degr\'{e} $3$ \`{a} coefficients r\'{e}els avec $\delta>0$. Ce polyn\^{o}me a
ses racines dans $\left\{  \operatorname{Re}z<\frac{1}{2}\right\}  $ si et
seulement si le polyn\^{o}me $Q(z)=P\left(  \frac{1}{2}+z\right)  $ est
stable. On a
\[
Q\left(  z\right)  =a_{0}z^{3}+a_{1}z^{2}+a_{2}z+a_{3},
\]
avec $a_{3}=P\left(  \frac{1}{2}\right)  $, $a_{2}=P^{\prime}\left(  \frac
{1}{2}\right)  $ et
\begin{align*}
\Delta_{2}  &  =%
\begin{vmatrix}
\beta+\gamma+\frac{3}{4}\delta & \delta\\
\alpha+\frac{\beta}{2}+\frac{\gamma}{4}+\frac{\delta}{8} & \gamma+\frac{3}%
{2}\delta
\end{vmatrix}
\\
&  =%
\begin{vmatrix}
\beta+\gamma+\frac{3}{4}\delta & \delta\\
\alpha-\frac{\gamma}{4}-\frac{\delta}{4} & \gamma+\delta
\end{vmatrix}
=%
\begin{vmatrix}
\beta+\gamma+\delta & \delta\\
\alpha & \gamma+\delta
\end{vmatrix}
.
\end{align*}
On en d\'{e}duit

\begin{proposition}
\label{LocDegre3}Soit $P(z)=\alpha+\beta z+\gamma z^{2}+\delta z^{3}$ un
polyn\^{o}me de degr\'{e} $3$ \`{a} coefficients r\'{e}els avec $\delta>0$. Ce
polyn\^{o}me a toutes ses racines dans $\left\{  \operatorname{Re}z<\frac
{1}{2}\right\}  $ si et seulement si on a
\begin{align}
P  &  \left(  \textstyle\frac{1}{2}\right)  >0,\quad P^{\prime}\left(
\textstyle\frac{1}{2}\right)  >0,\label{CritDeg3-a}\\
\Delta_{2}  &  =\left(  \gamma+\delta\right)  \left(  \beta+\gamma
+\delta\right)  -\alpha\delta>0. \label{CritDeg3-b}%
\end{align}

\end{proposition}

\subsection{Polyn\^{o}mes de degr\'{e} $4$}

Soit
\[
Q\left(  z\right)  =a_{0}z^{4}+a_{1}z^{3}+a_{2}z^{2}+a_{3}z+a_{4}%
\]
un polyn\^{o}me \`{a} coefficients r\'{e}els avec $a_{0}>0$. Les
d\'{e}terminants de Hurwitz sont
\begin{align*}
\Delta_{1}  &  =a_{1},\quad\Delta_{2}=
\begin{vmatrix}
a_{1} & a_{3}\\
a_{0} & a_{2}%
\end{vmatrix}
,\quad\Delta_{3}=
\begin{vmatrix}
a_{1} & a_{3} & 0\\
a_{0} & a_{2} & a_{4}\\
0 & a_{1} & a_{3}%
\end{vmatrix}
,\\
\Delta_{4}  &  =
\begin{vmatrix}
a_{1} & a_{3} & 0 & 0\\
a_{0} & a_{2} & a_{4} & 0\\
0 & a_{1} & a_{3} & 0\\
0 & 0 & a_{2} & a_{4}%
\end{vmatrix}
=a_{4}\Delta_{3}.
\end{align*}
Le polyn\^{o}me $Q$ est stable si et seulement tous les $\Delta_{j}$ sont
positifs. Les conditions $\Delta_{3}>0$ et $\Delta_{4}>0$ sont
\'{e}quivalentes \`{a}
\begin{equation}
\Delta_{3}>0,\quad a_{4}>0. \label{RouthDeg4-a}%
\end{equation}
D'autre part, on a
\[
\Delta_{3}=a_{3}\Delta_{2}-a_{1}^{2}a_{4}%
\]
et les conditions (\ref{RouthDeg4-a}) et $a_{3}>0$ entra\^{\i}nent donc
$\Delta_{2}>0$. Finalement, $\Delta_{2}=a_{1}a_{2}-a_{0}a_{3}$ et $a_{2}>0$
entra\^{\i}ne $a_{1}>0$ si les conditions pr\'{e}c\'{e}dentes sont
r\'{e}alis\'{e}es. Le crit\`{e}re de Routh--Hurwitz est ainsi \'{e}quivalent
au crit\`{e}re de Li\'{e}nard et Chipart, qui s'\'{e}crit ici

\begin{proposition}
\label{CritLC4}Le polyn\^{o}me \`{a} coefficients r\'{e}els%
\[
Q\left(  z\right)  =a_{0}z^{4}+a_{1}z^{3}+a_{2}z^{2}+a_{3}z+a_{4}%
\]
($a_{0}>0$) est stable si et seulement si on a
\begin{equation}
a_{2}>0,\quad a_{3}>0,\quad a_{4}>0,\quad\Delta_{3}>0. \label{LCDeg4}%
\end{equation}

\end{proposition}

Soit
\[
P(z)=\alpha+\beta z+\gamma z^{2}+\delta z^{3}+\varepsilon z^{4}%
\]
un polyn\^{o}me de degr\'{e} $4$ \`{a} coefficients r\'{e}els avec
$\varepsilon>0$. Soit
\[
Q\left(  z\right)  =a_{0}z^{4}+a_{1}z^{3}+a_{2}z^{2}+a_{3}z+a_{4}%
\]
le polyn\^{o}me $Q(z)=P\left(  \frac{1}{2}+z\right)  $. On a
\begin{align*}
a_{4}  &  =P\left(  \textstyle\frac{1}{2}\right)  =\textstyle\alpha+\frac
{1}{2}\beta+\frac{1}{4}\gamma+\frac{1}{8}\delta+\frac{1}{16}\varepsilon,\\
a_{3}  &  =P^{\prime}\left(  \textstyle\frac{1}{2}\right)  =\textstyle\beta
+\gamma+\frac{3}{4}\delta+\frac{1}{2}\varepsilon,\\
a_{2}  &  =\textstyle\frac{1}{2}P^{\prime\prime}\left(  \frac{1}{2}\right)
=\textstyle\gamma+\frac{3}{2}\delta+\frac{3}{2}\varepsilon,\\
a_{1}  &  =\textstyle\frac{1}{6}P^{\prime\prime\prime}\left(  \frac{1}%
{2}\right)  =\delta+2\varepsilon,\\
a_{0}  &  =\varepsilon.
\end{align*}
Le polyn\^{o}me de Hurwitz $\Delta_{3}$ relatif \`{a} $Q$ est
\begin{align*}
\Delta_{3}  &  =%
\begin{vmatrix}
a_{1} & a_{3} & 0\\
a_{0} & a_{2} & a_{4}\\
0 & a_{1} & a_{3}%
\end{vmatrix}
=a_{1}a_{2}a_{3}-a_{1}^{2}a_{4}-a_{0}a_{3}^{2}\\
&  =\textstyle\left(  \delta+2\varepsilon\right)  \left(  \gamma+\frac{3}%
{2}\delta+\frac{3}{2}\varepsilon\right)  \left(  \beta+\gamma+\frac{3}%
{4}\delta+\frac{1}{2}\varepsilon\right) \\
&  \qquad-\left(  \delta+2\varepsilon\right)  ^{2}\left(  \textstyle\alpha
+\frac{\beta}{2}+\frac{\gamma}{4}+\frac{\delta}{8}+\frac{\varepsilon}%
{16}\right)  -\varepsilon\left(  \textstyle\beta+\gamma+\frac{3}{4}%
\delta+\frac{1}{2}\varepsilon\right)  ^{2}.
\end{align*}
On a
\begin{equation}
\Delta_{3}=\left(  \varepsilon+\delta+\gamma+\beta\right)  \left[  \left(
\varepsilon+\delta+\gamma\right)  \left(  \varepsilon+\delta\right)
-\varepsilon\beta\right]  -\left(  2\varepsilon+\delta\right)  ^{2}\alpha.
\label{CRH5}%
\end{equation}

On a finalement, en appliquant cette relation et la proposition \ref{CritLC4}:

\begin{proposition}
\label{LocDegre4}Soit
\[
P(z)=\alpha+\beta z+\gamma z^{2}+\delta z^{3}+\varepsilon z^{4}%
\]
un polyn\^{o}me de degr\'{e} $4$ \`{a} coefficients r\'{e}els avec
$\varepsilon>0$. Les racines du polyn\^{o}me $P$ sont toutes situ\'{e}es dans
$\left\{  \operatorname{Re}z<\frac{1}{2}\right\}  $ si et seulement si $P$
v\'{e}rifie les conditions%
\begin{align}
P  &  \left(  \textstyle\frac{1}{2}\right)  >0,\quad P^{\prime}\left(
\textstyle\frac{1}{2}\right)  >0,\quad P^{\prime\prime}\left(  \textstyle\frac
{1}{2}\right)  >0,\label{CritDeg4-a}\\
\Delta_{3}  &  \equiv\left(  \varepsilon+\delta+\gamma+\beta\right)  \left[
\left(  \varepsilon+\delta+\gamma\right)  \left(  \varepsilon+\delta\right)
-\varepsilon\beta\right]  -\left(  2\varepsilon+\delta\right)  ^{2}\alpha>0.
\label{CritDeg4-b}%
\end{align}

\end{proposition}

\section{\label{Tables}Tables}

\numberwithin{equation}{subsection}

On trouvera ci-dessous pour les types indiqu\'{e}s de domaines born\'{e}s
sym\'{e}\-tri\-ques :

\begin{itemize}
\item le polyn\^{o}me de Hua
\[
\chi(s)=\underset{j=1}{\overset{r}{\prod}}\left(  \textstyle
s+1+\left(  j-1\right)  \frac{a}{2}\right)  _{1+b+\left(  r-j\right)  a}\ ;
\]

\item les coefficients $C_{j}(\mu)$ de la d\'{e}composition%
\[
\chi\left(  k\mu\right)  =\sum_{j=0}^{d}\mu^{j}C_{d-j}\left(  \mu\right)
\left(  k+1\right)  _{j}\ ;
\]

\item le polyn\^{o}me repr\'{e}sentatif du noyau de Bergman de $\widehat
{\Omega}_{m}\left(  \mu\right)  $%
\[
P_{\mu}^{m}(\eta)=\sum_{j=0}^{d}(m+1)_{j}C_{d-j}(\mu)\mu^{j}\eta^{j}\ ;
\]

\item ses polyn\^{o}mes d\'{e}riv\'{e}s $\frac{\operatorname*{d}^{k}%
}{\operatorname*{d}\eta^{k}}P_{\mu}^{m}(\eta)$ ($1\leq k<d$) ;

\item leurs valeurs pour $\eta=\frac{1}{2}$ ;

\item en dimension $3$, le d\'{e}terminant de Hurwitz
\[
\Delta_{2}=\left(  \gamma+\delta\right)  \left(  \beta+\gamma+\delta\right)
-\alpha\delta
\]
associ\'{e} au polyn\^{o}me $P_{\mu}^{m}\left(  \frac{1}{2}+\eta\right)
=\alpha+\beta\eta+\gamma\eta^{2}+\delta\eta^{3}\ ;$

\item en dimension $4$, le d\'{e}terminant de Hurwitz
\[
\Delta_{3}\equiv\left(  \varepsilon+\delta+\gamma+\beta\right)  \left[
\left(  \varepsilon+\delta+\gamma\right)  \left(  \varepsilon+\delta\right)
-\varepsilon\beta\right]  -\left(  2\varepsilon+\delta\right)  ^{2}\alpha
\]
associ\'{e} au polyn\^{o}me $P_{\mu}^{m}\left(  \frac{1}{2}+\eta\right)
=\alpha+\beta\eta+\gamma\eta^{2}+\delta\eta^{3}+\varepsilon\eta^{4}\ ;$

\item \'{e}ventuellement, les valeurs ci-dessus pour des valeurs
particuli\`{e}res de $m$.
\end{itemize}

\`{A} partir de la dimension $3$, les calculs ont \'{e}t\'{e} faits \`{a}
l'aide de \textsc{Mathematica} ; ils peuvent \^{e}tre v\'{e}rifi\'{e}s avec
tout autre logiciel de calcul symbolique.

\subsection{Type $I_{1,1}$\label{TableI-1-1}}

Le domaine $\Omega$ est le disque unit\'{e} de $\mathbb{C}$. On a
\begin{align*}
\chi\left(  s\right)   &  =s+1,\\
C_{1}(\mu)  &  =1-\mu,\quad C_{0}(\mu)=1,\\
P_{\mu}^{m}(\eta)  &  =1-\mu+(m+1)\mu\eta,\\
\textstyle P_{\mu}^{m}\left(  \frac{1}{2}\right)   &  =1+\textstyle\frac
{m-1}{2}\mu.
\end{align*}

\subsection{Type $I_{1,2}$\label{TableI-1-2}}

Le domaine $\Omega$ est la boule unit\'{e} de $\mathbb{C}^{2}$. On a
\begin{align*}
\chi\left(  s\right)   &  =(s+1)(s+2),\\
C_{0}(\mu)  &  =1,\quad C_{1}(\mu)=3(1-\mu),\quad C_{2}(\mu)=(1-\mu)(2-\mu),\\
P_{\mu}^{m}(\eta)  &  =(1-\mu)(2-\mu)+3(m+1)\mu(1-\mu)\eta+(m+1)_{2}\mu
^{2}\eta^{2}.
\end{align*}

\subsubsection{Signe de $P_{\mu}^{m}\left(  \frac{1}{2}\right)  $}

On a
\[
q_{m}(\mu)=\textstyle P_{\mu}^{m}\left(  \frac{1}{2}\right)
=\textstyle(2+m\mu)\left(  1+\frac{m-3}{4}\mu\right)  .
\]

\begin{proposition}
\label{Signe-I-1-2(0)}Pour $m\geq3$, $P_{\mu}^{m}(\frac{1}{2})$ est positif
pour tout $\mu>0$. Pour $m\leq2$, le polyn\^{o}me $q_{m}$ admet une seule
racine positive $\mu_{m}$ ($\mu_{1}=2<\mu_{2}=4$) et $P_{\mu}^{m}(\frac{1}%
{2})$ est positif ou nul si et seulement si $0<\mu\leq\mu_{m}$.
\end{proposition}

\subsubsection{Signe de $\frac{\operatorname*{d}P_{\mu}^{m}}{\operatorname*{d}%
\eta}\left(  \frac{1}{2}\right)  $}%

\[
\textstyle\frac{1}{(m+1)\mu}\frac{\operatorname*{d}P_{\mu}^{m}}%
{\operatorname*{d}\eta}\left(  \frac{1}{2}\right)  =3+(m-1)\mu.
\]

\begin{proposition}
\label{Signe-I-1-2(1)}Pour tout $m\geq1$ et pour tout $\mu>0$, on a
$\frac{\operatorname*{d}P_{\mu}^{m}}{\operatorname*{d}\eta}\left(  \frac{1}%
{2}\right)  >0$.
\end{proposition}

\subsection{Type $I_{1,3}$\label{TableI-1-3}}

Le domaine $\Omega$ est la boule unit\'{e} de $\mathbb{C}^{3}$. On a
\begin{align*}
\chi\left(  s\right)   &  =(s+1)(s+2)\left(  s+3\right)  ,\\
C_{0}(\mu)  &  =1,\quad C_{1}(\mu)=6(1-\mu),\\
C_{2}(\mu)  &  =(1-\mu)(11-7\mu),\quad C_{3}\left(  \mu\right)  =\left(
1-\mu\right)  \left(  2-\mu\right)  \left(  3-\mu\right)  ,\\
P_{\mu}^{m}(\eta)  &  =(1-\mu)(2-\mu)\left(  3-\mu\right)  +(m+1)(1-\mu
)\left(  11-7\mu\right)  \mu\eta\\
&  \qquad+6(m+1)_{2}\left(  1-\mu\right)  \mu^{2}\eta^{2}+\left(  m+1\right)
_{3}\mu^{3}\eta^{3}.
\end{align*}

\subsubsection{Signe de $P_{\mu}^{m}\left(  \frac{1}{2}\right)  $}

On a
\begin{align*}
\textstyle P_{\mu}^{m}\left(  \frac{1}{2}\right)   &  =q_{m}(\mu)=\left[
4+(m-1)\mu\right]  r_{m}(\mu),\\
r_{m}(\mu)  &  =\textstyle\frac{1}{8}\left(  12+8(m-1)\mu+(m^{2}-5m-2)\mu
^{2}\right)  .
\end{align*}
\emph{Cas particuliers} ($1\leq m\leq5$) :
\begin{align*}
r_{1}(\mu)  &  =\textstyle\frac{3}{4}\left(  2-\mu^{2}\right)  ,\\
r_{2}(\mu)  &  =\textstyle\frac{1}{2}(3+2\mu-2\mu^{2}),\\
r_{3}(\mu)  &  =\textstyle\frac{1}{2}(3+4\mu-2\mu^{2}),\\
r_{4}(\mu)  &  =\textstyle\frac{3}{4}(2+4\mu-\mu^{2}),\\
r_{5}(\mu)  &  =\textstyle\frac{1}{4}\left(  6+16\mu-\mu^{2}\right)  .
\end{align*}
\emph{Racines positives de }$q_{m}$ :
\[
\mu_{1}=\sqrt{2}<\mu_{2}=\textstyle\frac{1+\sqrt{7}}{2}<\mu_{3}=1+\sqrt
{\frac{5}{2}}<\mu_{4}=2+\sqrt{6}<\mu_{5}=8+\sqrt{70}.
\]

\begin{proposition}
\label{Signe-I-1-3(0)}Pour $m\geq6$, $P_{\mu}^{m}(\frac{1}{2})$ est positif
pour tout $\mu>0$. Pour $m\leq5$, le polyn\^{o}me $q_{m}$ admet une seule
racine positive $\mu_{m}$ et $P_{\mu}^{m}(\frac{1}{2})$ est positif ou nul si
et seulement si $0<\mu\leq\mu_{m}$.
\end{proposition}

\subsubsection{Signe de $\frac{\operatorname*{d}P_{\mu}^{m}}{\operatorname*{d}%
\eta}\left(  \frac{1}{2}\right)  $}

On a
\[
q_{m}^{1}(\mu)=\textstyle\frac{1}{(m+1)\mu}\frac{\operatorname*{d}P_{\mu}^{m}%
}{\operatorname*{d}\eta}\left(  \frac{1}{2}\right)  =11+6\left(  m-1\right)
\mu+\frac{1}{4}\left(  3m^{2}-9m-2\right)  \mu^{2}.
\]
\emph{Cas particuliers} ($1\leq m\leq3$) :
\begin{align*}
q_{1}^{1}(\mu)  &  =\textstyle11-2\mu^{2},\\
q_{2}^{1}(\mu)  &  =\textstyle11+6\mu-2\mu^{2},\\
q_{3}^{1}(\mu)  &  =\textstyle11+12\mu-\frac{1}{2}\mu^{2}.
\end{align*}
\emph{Racines positives de} $q_{m}^{1}$ ($1\leq m\leq3$) :%
\begin{align*}
0  &  <\mu_{1}^{1}=\textstyle\sqrt{\frac{11}{2}}<\mu_{2}^{1}=\frac{3+\sqrt
{31}}{2}<\mu_{3}^{1}=12+\sqrt{166},\\
0  &  <\mu_{m}<\mu_{m}^{1}\qquad(1\leq m\leq3).
\end{align*}

\begin{proposition}
\label{Signe-I-1-3(1)}Pour $m\geq4$, on a $\frac{\operatorname*{d}P_{\mu}^{m}%
}{\operatorname*{d}\eta}\left(  \frac{1}{2}\right)  >0$ pour tout $\mu>0$.
Pour $m\leq3$, le polyn\^{o}me $q_{m}^{1}$ poss\`{e}de une seule racine
positive $\mu_{m}^{1}$ et est positif sur $\left[  0,\mu_{m}^{1}\right]  $.
\end{proposition}

\subsubsection{Signe de $\frac{\operatorname*{d}^{2}P_{\mu}^{m}}%
{\operatorname*{d}\eta^{2}}\left(  \frac{1}{2}\right)  $}%

\[
\textstyle\frac{1}{\left(  m+1\right)  _{2}\mu^{2}}\frac{\operatorname*{d}%
^{2}P_{\mu}^{m}}{\operatorname*{d}\eta^{2}}\left(  \frac{1}{2}\right)
=3\left(  4+(m-1)\mu\right)  .
\]

\begin{proposition}
\label{Signe-I-1-3(2)}On a
\[
\textstyle\frac{\operatorname*{d}^{2}P_{\mu}^{m}}{\operatorname*{d}\eta^{2}%
}\left(  \frac{1}{2}\right)  >0
\]
pour tout $\mu>0$ et tout $m\geq1$.
\end{proposition}

\subsubsection{D\'{e}terminant de Hurwitz $\Delta_{2}$}

Soient $P_{\mu}^{m}(\eta)=\alpha+\beta\eta+\gamma\eta^{2}+\delta\eta^{3}$ et
$R_{m}(\mu)$ le polyn\^{o}me d\'{e}fini par
\[
R_{m}(\mu)=\Delta_{2}=\left(  \gamma+\delta\right)  \left(  \beta
+\gamma+\delta\right)  -\alpha\delta.
\]
On a
\begin{align*}
S_{m}(\mu)  &  =\textstyle\frac{1}{(m+1)_{2}\mu^{3}}R_{m}(\mu)=\left[
4+(m-1)\mu\right]  (3+m\mu)T_{m}(\mu),\\
T_{m}(\mu)  &  =5m+4+(m^{2}-2m-2)\mu.
\end{align*}
\emph{Cas particuliers} ($1\leq m\leq2$) : $T_{1}(\mu)=9-3\mu$, $T_{2}%
(\mu)=14-2\mu$.

\begin{proposition}
\label{Signe-I-1-3(Q)}Pour $m\geq3$, on a $R_{m}(\mu)>0$ pour tout $\mu>0$.
Pour $m\leq2$, le polyn\^{o}me $R_{m}$ admet une unique racine positive
$\nu_{m}$ ($\nu_{1}=3$, $\nu_{2}=7$) et on a $\nu_{m}>\mu_{m}$.
\end{proposition}

\subsubsection{Localisation des racines de $P_{\mu}^{m}$}

\begin{theorem}
\label{Racines-I-1-3}Soit $\Omega$ la boule hermitienne de dimension $3$ et
soit $P_{\mu}^{m}$ le polyn\^{o}me repr\'{e}sentatif du noyau de Bergman de
$\widehat{\Omega}_{m}\left(  \mu\right)  $. Pour $m\geq6$, les racines du
polyn\^{o}me $P_{\mu}^{m}$ sont dans le demi-plan $\left\{  \operatorname{Re}%
\eta\leq\frac{1}{2}\right\}  $ quel que soit $\mu>0$. Pour $m\leq5$, les
racines du polyn\^{o}me $P_{\mu}^{m}$ sont toutes dans le demi-plan $\left\{
\operatorname{Re}\eta\leq\frac{1}{2}\right\}  $ si et seulement si $0<\mu
\leq\mu_{m}$, o\`{u} $\mu_{m}$ d\'{e}signe l'unique racine positive de
$q_{m}(\mu)=P_{\mu}^{m}\left(  \frac{1}{2}\right)  $.
\end{theorem}

\begin{proof}
En rassemblant les r\'{e}sultats des propositions \ref{Signe-I-1-3(0)},
\ref{Signe-I-1-3(1)}, \ref{Signe-I-1-3(Q)} et en appliquant la proposition
\ref{LocDegre3}, on conclut que toutes les racines du polyn\^{o}me $P_{\mu
}^{m}$ sont dans le demi-plan $\left\{  \operatorname{Re}\eta<\frac{1}%
{2}\right\}  $ si et seulement si $0<\mu<\mu_{m}$. Comme l'ensemble des
racines varie contin\^{u}ment en fonction de $\mu$, les racines du
polyn\^{o}me $P_{\mu}^{m}$ sont dans le demi-plan $\left\{  \operatorname{Re}%
\eta\leq\frac{1}{2}\right\}  $ si et seulement si $0<\mu\leq\mu_{m}$.
\end{proof}

\subsection{Type $IV_{3}$\label{TableIV-3}}

Le domaine $\Omega$ est la boule de Lie de dimension $3$. Les invariants
num\'{e}riques sont $a=1$, $b=0$, $r=2$. Le polyn\^{o}me de Hua est%
\[
\textstyle\chi\left(  s\right)  =\left(  s+1\right)  \left(  s+\frac{3}%
{2}\right)  \left(  s+2\right)  .
\]
Les coefficients de la d\'{e}composition de $\chi(k\mu)$ sont
\begin{align*}
C_{0}(\mu)  &  =1,\quad C_{1}(\mu)=\textstyle3\left(  \frac{3}{2}-2\mu\right)
,\\
C_{2}(\mu)  &  =\textstyle(1-\mu)\left(  \frac{13}{2}-7\mu\right)  ,\\
C_{3}\left(  \mu\right)   &  =\textstyle\left(  1-\mu\right)  \left(
2-\mu\right)  \left(  \frac{3}{2}-\mu\right)  .
\end{align*}
Le polyn\^{o}me repr\'{e}sentatif du noyau de Bergman de $\widehat{\Omega}%
_{m}\left(  \mu\right)  $ est
\begin{align*}
P_{\mu}^{m} (\eta)=(  &  \textstyle1-\mu)(2-\mu)\left(  \frac{3}{2}%
-\mu\right)  +(m+1)\mu(1-\mu)\left(  \frac{13}{2}-7\mu\right)  \eta\\
&  +\textstyle3(m+1)_{2}\left(  \frac{3}{2}-2\mu\right)  \mu^{2}\eta
^{2}+\left(  m+1\right)  _{3}\mu^{3}\eta^{3}.
\end{align*}

\subsubsection{Signe de $P_{\mu}^{m}(\frac{1}{2})$}

On a
\begin{align*}
\textstyle P_{\mu}^{m}\left(  \frac{1}{2}\right)   &  =q_{m}(\mu)=\left[
3+(m-1)\mu\right]  r_{m}(\mu),\\
r_{m}(\mu)  &  =\textstyle\frac{1}{8}\left(  8+6(m-1)\mu+(m^{2}-5m-2)\mu
^{2}\right)  .
\end{align*}
\emph{Cas particuliers} ($1\leq m\leq5$) :%
\begin{align*}
r_{1}(\mu)  &  =\textstyle\frac{1}{4}\left(  4-3\mu^{2}\right)  ,\\
r_{2}(\mu)  &  =\textstyle\frac{1}{4}(4+3\mu-4\mu^{2}),\\
r_{3}(\mu)  &  =\textstyle\frac{1}{2}(1+2\mu)(2-\mu),\\
r_{4}(\mu)  &  =\textstyle\frac{1}{4}\left(  4+9\mu-3\mu^{2}\right)  ,\\
r_{5}(\mu)  &  =\textstyle\frac{1}{4}(4+12\mu-\mu^{2}).
\end{align*}
\emph{Racines positives de} $q_{m}$ ($1\leq m\leq5$) :%
\[
\mu_{1}=\textstyle\frac{2}{\sqrt{3}}<\mu_{2}=\frac{3+\sqrt{73}}{8}<\mu
_{3}=2<\mu_{4}=\frac{9+\sqrt{129}}{6}<\mu_{5}=2(3+\sqrt{10}).
\]

\begin{proposition}
\label{Signe-IV-3(0)}Pour $m\geq6$, $P_{\mu}^{m}(\frac{1}{2})$ est positif
pour tout $\mu>0$. Pour $m\leq5$, le polyn\^{o}me $q_{m}$ admet une seule
racine positive $\mu_{m}$ et $P_{\mu}^{m}(\frac{1}{2})$ est positif ou nul si
et seulement si $0<\mu\leq\mu_{m}$.
\end{proposition}

\subsubsection{Signe de $\frac{\operatorname*{d}P_{\mu}^{m}}{\operatorname*{d}%
\eta}\left(  \frac{1}{2}\right)  $}

On a
\[
q_{m}^{1}(\mu)=\textstyle\frac{1}{(m+1)\mu}\frac{\operatorname*{d}P_{\mu}^{m}%
}{\operatorname*{d}\eta}\left(  \frac{1}{2}\right)  =\frac{13}{2}+\frac{9}%
{2}(m-1)\mu+\frac{1}{4}\left(  3m^{2}-9m-2\right)  \mu^{2}.
\]
\emph{Cas particuliers} ($1\leq m\leq3$) :%
\begin{align*}
q_{1}^{1}(\mu)  &  =\textstyle\frac{13}{2}-2\mu^{2},\\
q_{2}^{1}(\mu)  &  =\textstyle\frac{1}{2}(1+\mu)(13-4\mu),\\
q_{3}^{1}(\mu)  &  =\textstyle\frac{13}{2}+9\mu-\frac{1}{2}\mu^{2}.
\end{align*}
\emph{Racines positives de} $q_{m}^{1}$ ($1\leq m\leq3$) :%
\begin{align*}
0  &  <\mu_{1}^{1}=\textstyle\frac{\sqrt{13}}{2}<\mu_{2}^{1}=\frac{13}{4}%
<\mu_{3}^{1}=9+\sqrt{94},\\
0  &  <\mu_{m}<\mu_{m}^{1}\qquad(1\leq m\leq3).
\end{align*}

\begin{proposition}
\label{Signe-IV-3(1)}Pour $m\geq4$, on a $\frac{\operatorname*{d}P_{\mu}^{m}%
}{\operatorname*{d}\eta}\left(  \frac{1}{2}\right)  >0$ pour tout $\mu\geq0$.
Pour $m\leq3$, le polyn\^{o}me $q_{m}^{1}$ poss\`{e}de une seule racine
positive $\mu_{m}^{1}$ et est positif sur $\left[  0,\mu_{m}^{1}\right]  $.
\end{proposition}

\subsubsection{Signe de $\frac{\operatorname*{d}^{2}P_{\mu}^{m}}%
{\operatorname*{d}\eta^{2}}\left(  \frac{1}{2}\right)  $}%

\[
\textstyle\frac{1}{\left(  m+1\right)  _{2}\mu^{2}}\frac{\operatorname*{d}%
^{2}P_{\mu}^{m}}{\operatorname*{d}\eta^{2}}\left(  \frac{1}{2}\right)
=3(3+(m-1)\mu).
\]

\begin{proposition}
\label{Signe-IV-3(2)}On a $\frac{\operatorname*{d}^{2}P_{\mu}^{m}%
}{\operatorname*{d}\eta^{2}}\left(  \frac{1}{2}\right)  >0$ pour tout $m\geq1$
et tout $\mu>0$.
\end{proposition}

\subsubsection{D\'{e}terminant de Hurwitz $\Delta_{2}$}

Soient $P_{\mu}^{m}(\eta)=\alpha+\beta\eta+\gamma\eta^{2}+\delta\eta^{3}$ et
$R_{m}(\mu)$ le polyn\^{o}me d\'{e}fini par
\[
R_{m}(\mu)=\Delta_{2}=\left(  \gamma+\delta\right)  \left(  \beta
+\gamma+\delta\right)  -\alpha\delta.
\]
On a
\begin{align*}
S_{m}(\mu)  &  =\textstyle\frac{1}{(m+1)_{2}\mu^{3}}R_{m}(\mu)=\left[
3+(m-1)\mu\right]  s_{m}(\mu),\\
s_{m}(\mu)  &  =\textstyle\frac{35m+27}{4}+\frac{3(m-1)(4m+3)}{2}\mu+m\left(
m-1\right)  \left(  m^{2}-2m-2\right)  \mu^{2}.
\end{align*}
\emph{Cas particuliers} ($1\leq m\leq2$) :
\begin{align*}
s_{1}(\mu)  &  =\textstyle\frac{1}{2}(31-6\mu^{2}),\\
s_{2}(\mu)  &  =\textstyle\frac{1}{4}\left(  97+66\mu-16\mu^{2}\right)  .
\end{align*}

\begin{proposition}
\label{Signe-IV-3(Q)}Pour $m\geq3$, on a $R_{m}(\mu)>0$ pour tout $\mu>0$.
Pour $1\leq m\leq2$, le polyn\^{o}me $R_{m}$ admet une unique racine positive
$\nu_{m}$ et on a $\nu_{m}>\mu_{m}$.
\end{proposition}

\subsubsection{Localisation des racines de $P_{\mu}^{m}$}

\begin{theorem}
\label{Racines-IV-3}Soit $\Omega$ un domaine sym\'{e}trique de type $IV_{3}$
et soit $P_{\mu}^{m}$ le polyn\^{o}me repr\'{e}sentatif du noyau de Bergman de
$\widehat{\Omega}_{m}\left(  \mu\right)  $. Pour $m\geq6$, les racines du
polyn\^{o}me $P_{\mu}^{m}$ sont dans le demi-plan $\left\{  \operatorname{Re}%
\eta\leq\frac{1}{2}\right\}  $ quel que soit $\mu>0$. Pour $m\leq5$, les
racines du polyn\^{o}me $P_{\mu}^{m}$ sont toutes dans le demi-plan $\left\{
\operatorname{Re}\eta\leq\frac{1}{2}\right\}  $ si et seulement si $0<\mu
\leq\mu_{m},$o\`{u} $\mu_{m}$ d\'{e}signe l'unique racine positive de
$q_{m}(\mu)=P_{\mu}^{m}\left(  \frac{1}{2}\right)  $.
\end{theorem}

\begin{proof}
En rassemblant les r\'{e}sultats des propositions \ref{Signe-IV-3(0)},
\ref{Signe-IV-3(1)}, \ref{Signe-IV-3(Q)} et en appliquant la proposition
\ref{LocDegre3}, on voit que toutes les racines du polyn\^{o}me $P_{\mu}^{m}$
sont dans le demi-plan $\left\{  \operatorname{Re}\eta<\frac{1}{2}\right\}  $
si et seulement si $0<\mu<\mu_{m}$. Comme l'ensemble des racines varie
contin\^{u}ment en fonction de $\mu$, les racines du polyn\^{o}me $P_{\mu}%
^{m}$ sont dans le demi-plan $\left\{  \operatorname{Re}\eta\leq\frac{1}%
{2}\right\}  $ si et seulement si $0<\mu\leq\mu_{m}$.
\end{proof}

\subsection{Type $I_{1,4}$\label{TableI-1-4}}

Le domaine $\Omega$ est la boule unit\'{e} de $%
\mathbb{C}
^{4}$. On a
\begin{align*}
\chi\left(  s\right)   &  =(s+1)(s+2)\left(  s+3\right)  (s+4),\\
C_{0}(\mu)  &  =1,\quad C_{1}(\mu)=10(1-\mu),\quad C_{2}(\mu)=5(1-\mu
)(7-5\mu),\\
C_{3}(\mu)  &  =5(1-\mu)(5-3\mu)(2-\mu),\\
C_{4}(\mu)  &  =\left(  1-\mu\right)  \left(  2-\mu\right)  \left(
3-\mu\right)  (4-\mu),\\
P_{\mu}^{m}(\eta)  &  =(1-\mu)(2-\mu)\left(  3-\mu\right)  (4-\mu)\\
&  \quad+5(m+1)(1-\mu)(5-3\mu)(2-\mu)\mu\eta\\
&  \quad+5(m+1)_{2}\left(  1-\mu\right)  (7-5\mu)\mu^{2}\eta^{2}\\
&  \quad+10\left(  m+1\right)  _{3}(1-\mu)\mu^{3}\eta^{3}+\left(  m+1\right)
_{4}\mu^{4}\eta^{4}.
\end{align*}

\subsubsection{Signe de $P_{\mu}^{m}\left(  \frac{1}{2}\right)  $}

On a%
\begin{align*}
q_{m}(\mu)  &  =\textstyle P_{\mu}^{m}\left(  \frac{1}{2}\right)
=(4+m\mu)r_{m}(\mu),\\
r_{m}(\mu)  &  =\textstyle6+\frac{19m-25}{4}\mu+m(m-5)\mu^{2}+\frac
{m^{3}-10m^{2}+15m+10}{16}\mu^{3}.
\end{align*}
\emph{Cas particuliers} ($1\leq m\leq7$) :
\begin{align*}
r_{1}(\mu)  &  =\textstyle\left(  \mu-4\right)  \left(  \mu^{2}-\frac{3}%
{2}\right)  ,\\
r_{2}(\mu)  &  =\textstyle6+\frac{13}{4}\mu-6\mu^{2}+\frac{1}{2}\mu^{3},\\
r_{3}(\mu)  &  =\textstyle6+8\mu-6\mu^{2}-\frac{1}{2}\mu^{3},\\
r_{4}(\mu)  &  =\textstyle6+\frac{51}{4}\mu-4\mu^{2}-\frac{13}{8}\mu^{3},\\
r_{5}(\mu)  &  =\textstyle6+\frac{35}{2}\mu-\frac{5}{2}\mu^{3},\\
r_{6}(\mu)  &  =\textstyle6+\frac{89}{4}\mu+6\mu^{2}-\frac{11}{4}\mu^{3},\\
r_{7}(\mu)  &  =\textstyle6+27\mu+14\mu^{2}-2\mu^{3}.
\end{align*}
\emph{Racines positives des polyn\^{o}mes }$q_{m}$ ($1\leq m\leq7$)\medskip\ :%

\begin{tabular}
[c]{c|ccccccc}%
$m$ & $1$ & $2$ & $3$ & $4$ & $5$ & $6$ & $7$\\\hline
$\mu_{m}$ & $\sqrt{\frac{3}{2}}$ & $1.41518$ & $1.68819$ & $2.10335$ &
$2.8029$ & $4.22107$ & $8.60867$\\\hline
$\mu_{m,2}$ & $4$ & $11.333$ &  &  &  &  &
\end{tabular}

\begin{proposition}
\label{Signe-I-1-4(0)}Pour $m\geq8$, $P_{\mu}^{m}(\frac{1}{2})$ est positif
pour tout $\mu>0$. Pour $m\leq2$, le polyn\^{o}me $q_{m}$ admet deux racines
positives $\mu_{m}=\mu_{m,1}<\mu_{m,2}$, et $P_{\mu}^{m}(\frac{1}{2})$ est
positif ou nul si et seulement si $0<\mu\leq\mu_{m}\ $ou $\mu\geq\mu_{m,2}$.
Pour $3\leq m\leq7$, le polyn\^{o}me $q_{m}$ admet une seule racine positive
$\mu_{m}$ et $P_{\mu}^{m}(\frac{1}{2})$ est positif ou nul si et seulement si
$0<\mu\leq\mu_{m}$.
\end{proposition}

\subsubsection{Signe de $\frac{\operatorname*{d}P_{\mu}^{m}}{\operatorname*{d}%
\eta}\left(  \frac{1}{2}\right)  $}

On a
\begin{align*}
q_{m}^{1}(\mu)  &  =\textstyle\frac{1}{(m+1)\mu}\frac{\operatorname*{d}P_{\mu
}^{m}}{\operatorname*{d}\eta}\left(  \frac{1}{2}\right)  =(5+(m-1)\mu
)r_{m}^{1}(\mu),\\
r_{m}^{1}(\mu)  &  =\textstyle10+5(m-1)\mu+\frac{m^{2}-5m-4}{2}\mu^{2}.
\end{align*}
\emph{Cas particuliers} ($1\leq m\leq5$) :
\begin{align*}
r_{1}^{1}(\mu)  &  =2(5-2\mu^{2}),\\
r_{2}^{1}(\mu)  &  =5(1+\mu)(2-\mu),\\
r_{3}^{1}(\mu)  &  =5(2+2\mu-\mu^{2}),\\
r_{4}^{1}(\mu)  &  =10+15\mu-4\mu^{2},\\
r_{5}^{1}(\mu)  &  =2(5+10\mu-\mu^{2}).
\end{align*}
\emph{Racines de }$q_{m}^{1}$ ($1\leq m\leq5$) : \medskip%

\begin{tabular}
[c]{c|ccccc}%
$m$ & $1$ & $2$ & $3$ & $4$ & $5$\\\hline
$\mu_{m}^{1}$ & $\sqrt{\frac{5}{2}}$ & $2$ & $\sqrt{3}+1$ & $\frac{1}%
{8}\left(  \sqrt{385}+15\right)  $ & $\sqrt{30}+5$%
\end{tabular}
\medskip

\noindent On a $0<\mu_{m}<\mu_{m}^{1}$ ($1\leq m\leq5$).

\begin{proposition}
\label{Signe-I-1-4(1)}Pour $m\geq6$, on a $\frac{\operatorname*{d}P_{\mu}^{m}%
}{\operatorname*{d}\eta}\left(  \frac{1}{2}\right)  >0$ pour tout $\mu>0$.
Pour $m\leq5$, le polyn\^{o}me $q_{m}^{1}$ poss\`{e}de une seule racine
positive $\mu_{m}^{1}$ et est positif sur $\left[  0,\mu_{m}^{1}\right]  $.
\end{proposition}

\subsubsection{Signe de $\frac{\operatorname*{d}^{2}P_{\mu}^{m}}%
{\operatorname*{d}\eta^{2}}\left(  \frac{1}{2}\right)  $}

On a%
\[
q_{m}^{2}(\mu)=\textstyle\frac{1}{\left(  m+1\right)  _{2}\mu^{2}}%
\frac{\operatorname*{d}^{2}P_{\mu}^{m}}{\operatorname*{d}\eta^{2}}\left(
\frac{1}{2}\right)  =\left(  3m^{2}-9m-4\right)  \mu^{2}+30\left(  m-1\right)
\mu+70.
\]
\emph{Cas particuliers} ($1\leq m\leq3$) :%
\begin{align*}
q_{1}^{2}(\mu)  &  =\textstyle10\left(  7-\mu^{2}\right)  ,\\
q_{2}^{2}(\mu)  &  =\textstyle10\left(  7+3\mu-\mu^{2}\right)  ,\\
q_{3}^{2}(\mu)  &  =\textstyle2\left(  35+30\mu-2\mu^{2}\right)  .
\end{align*}
\emph{Racines positives de }$q_{m}^{2}$ ($1\leq m\leq3$) :\medskip%

\begin{tabular}
[c]{c|ccc}%
$m$ & $1$ & $2$ & $3$\\\hline
$\mu_{m}^{2}$ & $\sqrt{7}$ & $\frac{1}{2}\sqrt{37}+\frac{3}{2}$ & $\frac{1}%
{2}\sqrt{295}+\frac{15}{2}$%
\end{tabular}

\medskip\noindent On a $0<\mu_{m}<\mu_{m}^{1}<\mu_{m}^{2}$ ($1\leq m\leq3$) et
$\mu_{m}^{2}<\mu_{m,2}$ ($1\leq m\leq2$).

\begin{proposition}
\label{Signe-I-1-4(2)}Pour $m\geq4$, on a $\frac{\operatorname*{d}^{2}P_{\mu
}^{m}}{\operatorname*{d}\eta^{2}}\left(  \frac{1}{2}\right)  >0$ pour tout
$\mu>0$. Pour $m\leq3$, le polyn\^{o}me $q_{m}^{2}$ poss\`{e}de une seule
racine positive $\mu_{m}^{2}$ et est positif sur $\left[  0,\mu_{m}%
^{2}\right]  $.
\end{proposition}

\subsubsection{Signe de $\frac{\operatorname*{d}^{3}P_{\mu}^{m}}%
{\operatorname*{d}\eta^{3}}\left(  \frac{1}{2}\right)  $}%

\[
\textstyle\frac{1}{\left(  m+1\right)  _{3}\mu^{3}}\frac{\operatorname*{d}%
^{3}P_{\mu}^{m}}{\operatorname*{d}\eta^{3}}\left(  \frac{1}{2}\right)
=12\left(  \left(  m-1\right)  \mu+5\right)  .
\]

\begin{proposition}
\label{Signe-I-1-4(3)}On a $\frac{\operatorname*{d}^{3}P_{\mu}^{m}%
}{\operatorname*{d}\eta^{3}}\left(  \frac{1}{2}\right)  >0$ pour tout $\mu>0$
et tout $m\geq1$.
\end{proposition}

\subsubsection{D\'{e}terminant de Hurwitz $\Delta_{3}$}

Pour $P_{\mu}^{m}(\eta)=\alpha+\beta\eta+\gamma\eta^{2}+\delta\eta
^{3}+\varepsilon\eta^{4}$, soit
\[
F_{m}(\mu)=\Delta_{3}\equiv\left(  \varepsilon+\delta+\gamma+\beta\right)
\left[  \left(  \varepsilon+\delta+\gamma\right)  \left(  \varepsilon
+\delta\right)  -\varepsilon\beta\right]  -\left(  2\varepsilon+\delta\right)
^{2}\alpha.
\]
On a%
\begin{align*}
G_{m}(\mu)  &  =\textstyle\frac{1}{(m+1)(m+1)_{3}\mu^{6}}F_{m}(\mu)=\left(
5+\left(  m-1\right)  \mu\right)  ^{2}\left(  4+m\mu\right)  S_{m}\left(
\mu\right) \\
S_{m}(\mu)  &  =\sum_{j=0}^{3}s_{j}(m)\mu^{j},\\
s_{0}(m)  &  =2\left(  53+140m+63m^{2}\right)  ,\\
s_{1}(m)  &  =-75-189m+55m^{2}+81m^{3},\\
s_{2}(m)  &  =2m\left(  -5-52m-15m^{2}+8m^{3}\right)  ,\\
s_{3}(m)  &  =5+15m+20m^{2}-4m^{3}-5m^{4}+m^{5}.
\end{align*}
Pour $m>4$, les coefficients $s_{j}(m)$ sont tous positifs.

\noindent\emph{Cas particuliers :}%
\begin{align*}
S_{1}(\mu)  &  =32\left(  \mu-2\right)  \left(  \mu+2\right)  \left(
\mu-4\right)  ,\\
S_{2}(\mu)  &  =1170+415\mu-521\mu^{2}+35\mu^{3},\\
S_{3}(\mu)  &  =40\left(  4-\mu\right)  \left(  16\mu+\mu^{2}+13\right)  ,\\
S_{4}(\mu)  &  =3242+5233\mu+354\mu^{2}-127\mu^{3}.
\end{align*}
\noindent\emph{Racines positives des polyn\^{o}mes }$S_{m}$ ($1\leq m\leq4$) :
\medskip%

\begin{tabular}
[c]{c|cccc}%
$m$ & $1$ & $2$ & $3$ & $4$\\\hline
$\sigma_{m}=\sigma_{m,1}\simeq$ & $2$ & $2.15132$ & $4$ & $8.19532$\\\hline
$\sigma_{m,2}\simeq$ & $4$ & $13.8558$ &  &
\end{tabular}
\medskip

\noindent On a $\mu_{m}^{1}<\sigma_{m,1}<\sigma_{m,2}$ ($m=1,2$) et $\mu
_{m}^{1}<\sigma_{m}$ ($m=3,4$), o\`{u} $\mu_{m}^{1}$ est la racine positive de
$q_{m}^{1}$.

\begin{proposition}
\label{Hurwitz3-I-1-4}Pour $m>4$, on a $F_{m}(\mu)>0$ pour tout $\mu>0$. Pour
$m\leq4$, on a $F_{m}(\mu)>0$ et $\frac{\operatorname*{d}P_{\mu}^{1}%
}{\operatorname*{d}\eta}\left(  \frac{1}{2}\right)  \geq0$ si et seulement si
$\mu\leq\mu_{m}^{1}$.
\end{proposition}

\subsubsection{Localisation des racines de $P_{\mu}^{m}$}

\begin{theorem}
\label{LocTypeI-1-4}Les racines du polyn\^{o}me $P_{\mu}^{m}$ sont toutes
situ\'{e}es dans le demi-plan $\left\{  \operatorname{Re}\eta\leq\frac{1}%
{2}\right\}  $

\begin{itemize}
\item pour $m\geq8$ et pour tout $\mu>0$ ;

\item pour $1\leq m\leq7$, si et seulement si $0<\mu\leq\mu_{m}$, o\`{u}
$\mu_{m}$ est la plus petite racine positive de $q_{m}(\mu)=P_{\mu}^{m}\left(
\frac{1}{2}\right)  $.
\end{itemize}
\end{theorem}

\begin{proof}
En rassemblant les r\'{e}sultats des propositions (\ref{Signe-I-1-4(0)},
\ref{Signe-I-1-4(1)}, \ref{Signe-I-1-4(2)}) et en appliquant la proposition
\ref{Hurwitz3-I-1-4}, on voit que toutes les racines du polyn\^{o}me $P_{\mu
}^{m}$ sont dans le demi-plan $\left\{  \operatorname{Re}\eta<\frac{1}%
{2}\right\}  $ pour $0<\mu<\mu_{m}$. Comme l'ensemble des racines varie
contin\^{u}ment en fonction de $\mu$, les racines du polyn\^{o}me $P_{\mu}%
^{m}$ sont dans le demi-plan $\left\{  \operatorname{Re}\eta\leq\frac{1}%
{2}\right\}  $ si et seulement si $0<\mu\leq\mu_{m}$.
\end{proof}

\subsection{Type $IV_{4}$\label{TableIV-4}}

Le domaine $\Omega$ est la boule de Lie de dimension $4$. Les invariants
num\'{e}riques sont $a=2$, $b=0$, $r=2$. Le polyn\^{o}me de Hua est%
\[
\chi\left(  s\right)  =\left(  s+1\right)  \left(  s+2\right)  ^{2}(s+3).
\]
Les coefficients de la d\'{e}composition de $\chi(\mu k)$ sont
\begin{align*}
C_{0}(\mu)  &  =1,\quad C_{1}(\mu)=2(4-5\mu),\quad C_{2}(\mu)=(1-\mu
)(23-25\mu),\\
C_{3}\left(  \mu\right)   &  =(1-\mu)(7-5\mu)(4-3\mu),\quad C_{4}(\mu)=\left(
1-\mu\right)  \left(  2-\mu\right)  ^{2}(3-\mu).
\end{align*}
Le polyn\^{o}me repr\'{e}sentatif du noyau de Bergman de $\widehat{\Omega}%
_{m}\left(  \mu\right)  $ est
\begin{align*}
P_{\mu}^{m} (\eta)=(  &  1-\mu) \left(  2-\mu\right)  ^{2}(3-\mu
)+(m+1)(1-\mu)(7-5\mu)(4-3\mu)\mu\eta\\
&  +(m+1)_{2}(1-\mu)(23-25\mu)\mu^{2}\eta^{2}+2\left(  m+1\right)  _{3}%
(4-5\mu)\mu^{3}\eta^{3}\\
&  +\left(  m+1\right)  _{4}\mu^{4}\eta^{4}.
\end{align*}

\subsubsection{Signe de $P_{\mu}^{m}(\frac{1}{2})$}

Soit
\begin{align*}
q_{m}(\mu)  &  =P_{\mu}^{m}\left(  \textstyle\frac{1}{2}\right) \\
&  =12+\textstyle14(m-1)\mu+\frac{23}{4}m(m-3)\mu^{2}\\
&  \qquad+\textstyle(m-1)(m^{2}-5m-2)\mu^{3}+\frac{1}{16}\left(  m^{3}%
-10m^{2}+15m+10\right)  \mu^{4}.
\end{align*}
\emph{Cas particuliers} ($1\leq m\leq7$) :%
\begin{align*}
q_{1}(\mu)  &  =\textstyle\frac{1}{2}\left(  24-23\mu^{2}+2\mu^{4}\right)  ,\\
q_{2}(\mu)  &  =\textstyle\frac{1}{2}\left(  24+28\mu-23\mu^{2}-16\mu^{3}%
+2\mu^{4}\right)  ,\\
q_{3}(\mu)  &  =\textstyle\frac{1}{2}\left(  24+56\mu-32\mu^{3}-3\mu
^{4}\right)  ,\\
q_{4}(\mu)  &  =\textstyle\frac{1}{2}\left(  24+84\mu+46\mu^{2}-36\mu
^{3}-13\mu^{4}\right)  ,\\
q_{5}(\mu)  &  =\textstyle\frac{1}{2}\left(  24+112\mu+115\mu^{2}-16\mu
^{3}-25\mu^{4}\right)  ,\\
q_{6}(\mu)  &  =\textstyle\frac{1}{2}\left(  24+140\mu+207\mu^{2}+40\mu
^{3}-33\mu^{4}\right)  ,\\
q_{7}(\mu)  &  =\textstyle12+84\mu+161\mu^{2}+72\mu^{3}-14\mu^{4}.
\end{align*}
\emph{Racines positives des polyn\^{o}mes }$q_{m}$ ($1\leq m\leq7$) : \medskip

\noindent%
\begin{tabular}
[c]{c|ccccccc}%
$m$ & $1$ & $2$ & $3$ & $4$ & $5$ & $6$ & $7$\\\hline
$\mu_{m,1}$ & $1.07732$ & $1.21176$ & $1.41824$ & $1.74173$ & $2.29476$ &
$3.42405$ & $6.92986$\\\hline
$\mu_{m,2}$ & $3.21549$ & $9.08062$ &  &  &  &  &
\end{tabular}

\medskip\noindent On note $\mu_{m}=\mu_{m,1}$ ($1\leq m\leq7$).

\begin{proposition}
\label{Signe-IV-4(0)}Pour $m\geq8$, $P_{\mu}^{m}(\frac{1}{2})$ est positif
pour tout $\mu>0$. On a :

\begin{enumerate}
\item $P_{\mu}^{1}(\frac{1}{2})\geq0$ pour $0<\mu\leq\mu_{m,1}$ et $\mu
_{m,2}\leq\mu$ ($1\leq m\leq2$) ;

\item $P_{\mu}^{m}(\frac{1}{2})\geq0$ si $0<\mu\leq\mu_{m}=\mu_{m,1}$ ($3\leq
m\leq7$).
\end{enumerate}
\end{proposition}

\subsubsection{Signe de $\frac{\operatorname*{d}P_{\mu}^{m}}{\operatorname*{d}%
\eta}\left(  \frac{1}{2}\right)  $}

On a
\begin{align*}
q_{m}^{1}(\mu)  &  =\textstyle\frac{1}{(m+1)\mu}\frac{\operatorname*{d}P_{\mu
}^{m}}{\operatorname*{d}\eta}\left(  \frac{1}{2}\right)  =\left[
4+(m-1)\mu\right]  r_{m}^{1},\\
r_{m}^{1}(\mu)  &  =\textstyle7+4(m-1)\mu+\frac{1}{2}\left(  m^{2}%
-5m-4\right)  \mu^{2}.
\end{align*}
\emph{Cas particuliers} ($1\leq m\leq5$) :%
\begin{align*}
r_{1}^{1}(\mu)  &  =7-4\mu^{2},\\
r_{2}^{1}(\mu)  &  =7+4\mu-5\mu^{2},\\
r_{3}^{1}(\mu)  &  =7+8\mu-5\mu^{2},\\
r_{4}^{1}(\mu)  &  =(1+2\mu)(7-2\mu),\\
r_{5}^{1}(\mu)  &  = 7+16\mu-2\mu^{2} .
\end{align*}
\emph{Racines positives de }$q_{m}^{1}$ ($1\leq m\leq5$) : \medskip%

\begin{tabular}
[c]{c|ccccc}%
$m$ & $1$ & $2$ & $3$ & $4$ & $5$\\\hline
$\mu_{m}^{1}$ & $\frac{\sqrt{7}}{2}$ & $\simeq1.649$ & $\simeq2.22829$ &
$\frac{7}{2}$ & $\simeq8.41588$%
\end{tabular}

\medskip\noindent On a $0<\mu_{m,1}<\mu_{m}^{1}<\mu_{m,2}$ ($1\leq m\leq2$) et
$0<\mu_{m}<\mu_{m}^{1}$ ($3\leq m\leq5$).

\begin{proposition}
\label{Signe-IV-4(1)}Pour $m\geq6$, on a $\frac{\operatorname*{d}P_{\mu}^{m}%
}{\operatorname*{d}\eta}\left(  \frac{1}{2}\right)  >0$ pour tout $\mu\geq0$.
Pour $m\leq6$, le polyn\^{o}me $q_{m}^{1}$ poss\`{e}de une seule racine
positive $\mu_{m}^{1}$ et est positif sur $\left[  0,\mu_{m}^{1}\right]  $.
\end{proposition}

\subsubsection{Signe de $\frac{\operatorname*{d}^{2}P_{\mu}^{m}}%
{\operatorname*{d}\eta^{2}}\left(  \frac{1}{2}\right)  $}

On a
\[
q_{m}^{2}(\mu)=\textstyle\frac{1}{\left(  m+1\right)  _{2}\mu^{2}}%
\frac{\operatorname*{d}^{2}P_{\mu}^{m}}{\operatorname*{d}\eta^{2}}\left(
\frac{1}{2}\right)  =46+24(m-1)\mu+\left(  3m^{2}-9m-4\right)  \mu^{2}.
\]
\emph{Cas particuliers : }%
\begin{align*}
q_{1}^{2}(\mu)  &  =2(23-25\mu^{2}),\\
q_{2}^{2}(\mu)  &  =2\left(  23+12\mu-5\mu^{2}\right)  ,\\
q_{3}^{2}(\mu)  &  =2\left(  23+24\mu-2\mu^{2}\right)  .
\end{align*}
\emph{Racines positives de }$q_{m}^{2}$ ($1\leq m\leq3$) :\medskip%

\begin{tabular}
[c]{c|ccc}%
$m$ & $1$ & $2$ & $3$\\\hline
$\mu_{m}^{2}$ & $\simeq2.14476$ & $\simeq3.65764$ & $\simeq12.892$%
\end{tabular}

\medskip\noindent On a $\mu_{1}^{m}<\mu_{2}^{m}$ ($1\leq m\leq3$).

\begin{proposition}
\label{Signe-IV-4(2)}On a $\frac{\operatorname*{d}^{2}P_{\mu}^{m}%
}{\operatorname*{d}\eta^{2}}\left(  \frac{1}{2}\right)  >0$ pour tout $m\geq4$
et tout $\mu>0$. Pour $m\leq3$, le polyn\^{o}me $q_{m}^{2}$ poss\`{e}de une
seule racine positive $\mu_{m}^{2}$ et est positif si $0<\mu<\mu_{m}^{2}$. En
particulier, on a $\frac{\operatorname*{d}^{2}P_{\mu}^{m}}{\operatorname*{d}%
\eta^{2}}\left(  \frac{1}{2}\right)  >0$ pour $1\leq m\leq7$ et $0<\mu\leq
\mu_{m}$.
\end{proposition}

\subsubsection{Signe de $\frac{\operatorname*{d}^{3}P_{\mu}^{m}}%
{\operatorname*{d}\eta^{3}}\left(  \frac{1}{2}\right)  $}%

\[
\textstyle\frac{1}{\left(  m+1\right)  _{3}\mu^{3}}\frac{\operatorname*{d}%
^{3}P_{\mu}^{m}}{\operatorname*{d}\eta^{3}}\left(  \frac{1}{2}\right)
=12\left[  4+(m-1)\mu\right]  .
\]

\begin{proposition}
\label{Signe-IV-4(3)}On a $\frac{\operatorname*{d}^{3}P_{\mu}^{m}%
}{\operatorname*{d}\eta^{3}}\left(  \frac{1}{2}\right)  >0$ pour tout $m\geq1$
et tout $\mu>0$.
\end{proposition}

\subsubsection{D\'{e}terminant de Hurwitz $\Delta_{3}$}

Pour $P_{\mu}^{m}(\eta)=\alpha+\beta\eta+\gamma\eta^{2}+\delta\eta
^{3}+\varepsilon\eta^{4}$, soit
\[
F_{m}(\mu)=\Delta_{3}\equiv\left(  \varepsilon+\delta+\gamma+\beta\right)
\left[  \left(  \varepsilon+\delta+\gamma\right)  \left(  \varepsilon
+\delta\right)  -\varepsilon\beta\right]  -\left(  2\varepsilon+\delta\right)
^{2}\alpha.
\]
On a
\begin{align*}
G_{m}(\mu)  &  =\textstyle\frac{1}{(m+1)(m+1)_{3}\mu^{6}}F_{m}(\mu)={\left(
4+(m-1)\mu\right)  }^{2}S_{m}(\mu),\\
S_{m}(\mu)  &  =\sum_{j=0}^{4}s_{j}(m)\mu^{j},\\
s_{0}(m)  &  =160+481m+225m^{2},\\
s_{1}(m)  &  =16\left(  m-1\right)  \left(  9+31m+15m^{2}\right)  ,\\
s_{2}(m)  &  =2m\left(  -35-196m-22m^{2}+47m^{3}\right)  ,\\
s_{3}(m)  &  =8\left(  m-1\right)  \left(  -2-8m-15m^{2}-3m^{3}+2m^{4}\right)
,\\
s_{4}(m)  &  =m\left(  5+15m+20m^{2}-4m^{3}-5m^{4}+m^{5}\right)  .
\end{align*}
Pour $m>4$, les coefficients $s_{j}(m)$ sont tous positifs.

\noindent\emph{Cas particuliers : }%
\begin{align*}
S_{1}(\mu)  &  =2\left(  433-206{\mu}^{2}+16{\mu}^{4}\right)  ,\\
S_{2}(\mu)  &  =2\,\left(  1+\mu\right)  \left(  1011+37\mu-315{\mu}%
^{2}+35{\mu}^{3}\right)  ,\\
S_{3}(\mu)  &  =4\left(  907+1896\mu+672{\mu}^{2}-320{\mu}^{3}-30{\mu}%
^{4}\right)  ,\\
S_{4}(\mu)  &  =4\left(  1421+4476\mu+3674{\mu}^{2}+276{\mu}^{3}-127{\mu}%
^{4}\right)  .
\end{align*}
\emph{Racines positives des polyn\^{o}mes }$S_{m}$ ($1\leq m\leq4$) : \medskip%

\begin{tabular}
[c]{c|cccc}%
$m$ & $1$ & $2$ & $3$ & $4$\\\hline
$\sigma_{m}=\sigma_{m,1}\simeq$ & $1.62651$ & $2.12869$ & $3.22913$ &
$7.03204$\\\hline
$\sigma_{m,2}\simeq$ & $3.19835$ & $8.47286$ &  &
\end{tabular}
\medskip

\noindent On a $\mu_{m}^{1}<\sigma_{m,1}<\sigma_{m,2}$ ($m=1,2$) et $\mu
_{m}^{1}<\sigma_{m}$ ($m=3,4$), o\`{u} $\mu_{m}^{1}$ est la racine positive de
$q_{m}^{1}$.

\begin{proposition}
\label{Hurwitz3-IV-4}Pour $m>4$, on a $F_{m}(\mu)>0$ pour tout $\mu>0$. Pour
$m\leq4$, on a $F_{m}(\mu)>0$ et $\frac{\operatorname*{d}P_{\mu}^{1}%
}{\operatorname*{d}\eta}\left(  \frac{1}{2}\right)  \geq0$ si et seulement si
$\mu\leq\mu_{m}^{1}$.
\end{proposition}

\subsection{Localisation des racines de $P_{\mu}^{m}$}

\begin{theorem}
\label{LocTypeIV-4}Les racines du polyn\^{o}me $P_{\mu}^{m}$ sont toutes
situ\'{e}es dans le demi-plan $\left\{  \operatorname{Re}\eta\leq\frac{1}%
{2}\right\}  $

\begin{itemize}
\item pour $m\geq8$ et pour tout $\mu>0$ ;

\item pour $1\leq m\leq7$, si et seulement si $0<\mu\leq\mu_{m}$, o\`{u}
$\mu_{m}$ est la plus petite racine positive de $q_{m}(\mu)=P_{\mu}^{m}\left(
\frac{1}{2}\right)  $.
\end{itemize}
\end{theorem}

\begin{proof}
En rassemblant les r\'{e}sultats des propositions \ref{Signe-IV-4(0)},
\ref{Signe-IV-4(1)}, \ref{Signe-IV-4(2)}, \ref{Signe-IV-4(3)} et en appliquant
la proposition \ref{Hurwitz3-IV-4}, on voit que toutes les racines du
polyn\^{o}me $P_{\mu}^{m}$ sont dans le demi-plan $\left\{  \operatorname{Re}%
\eta<\frac{1}{2}\right\}  $ si et seulement si $0<\mu<\mu_{m}$. Comme
l'ensemble des racines varie contin\^{u}ment en fonction de $\mu$, les racines
du polyn\^{o}me $P_{\mu}^{m}$ sont dans le demi-plan $\left\{
\operatorname{Re}\eta\leq\frac{1}{2}\right\}  $ si et seulement si $0<\mu
\leq\mu_{m}$.
\end{proof}

\selectlanguage{french}

\end{document}